\documentclass[preprints,article,accept,moreauthors,pdftex]{mdpi} 

\firstpage{1} 
\pubvolume{xx}
\issuenum{1}
\articlenumber{5}
\pubyear{2019}
\copyrightyear{2019}
\externaleditor{Academic Editor: Georgy Sofronov}
\history{Received: 22nd of September 2020; Accepted: date; Published: date}

\usepackage{comment} 
\usepackage{amssymb} 
\graphicspath{{./Pictures/},{./Definitions/}}

\newcommand*{\doi}[1]{\href{http://dx.doi.org/#1}{doi: #1}}

\def\bP{\mathbf P}
\def\bE{\mathbf E}
\def\card{\text{card}}
\Title{Cross-entropy method in application to SIRC model}

\Author{Maria Katarzyna Stachowiak $^{1}$\orcidA{} and Krzysztof Józef Szajowski $^{2}$\orcidB{}}

\AuthorNames{Maria Katarzyna Stachowiak and Krzysztof Józef Szajowski}

\address{%
$^{1}$ \quad Faculty of Pure and Applied Mathematics, Wrocław University of Science and Technology, Wybrze\.ze Wyspia{\'n}skiego 27, 50-370 Wroc{\l}aw, Poland; mkwstachowiak@gmail.com\\
$^{2}$ \quad Faculty of Pure and Applied Mathematics, Wrocław University of Science and Technology, Wybrze\.ze Wyspia{\'n}skiego 27, 50-370 Wroc{\l}aw, Poland; Krzysztof.Szajowski@pwr.edu.pl}

\corres{Correspondence: mkwstachowiak@gmail.com(MKS); kszajowsk@gmail.com; Tel.: +48-71-3203185(KJS)}

\abstract{The study considers the usage of a probabilistic optimization method called Cross-Entropy (\textbf{CE}). This is the version of the Monte Carlo method created by Reuven Rubinstein (1997). It was developed in the context of determining rare events. Here we will present the way in which the \textbf{CE} method can be used for problems of optimization of epidemiological models, and more specifically the optimization of the SIRC (\textbf{S}usceptible - \textbf{I}nfectious - \textbf{R}ecovered - \textbf{C}ross-immune) model based on the functions supervising the care of specific groups in the model. With the help of weighted sampling, an attempt was made to find the fastest and most accurate version of the algorithm.}



\usepackage{booktabs}
    \usepackage{enumitem}
    \usepackage{lscape}
    \newlist{tableitems}{itemize}{1}
    \usepackage{mathabx}
    \setlist[tableitems]{nosep,
                         topsep=0pt,
                         partopsep=0pt,
                         leftmargin=1em,
                         label=$\sqbullet$
    }
\begin{document}


\section{\label{chapter:introduction}Introduction.} In this study our aim is to develop the possibilities of numerically solving variational problems that appear in epidemic dynamics models. The algorithm that we will use for this purpose will be a modified Monte Carlo method. In a few sentences, we recall those stages of the Monte Carlo method development and its improvement, which constitute the essence of the applied approach (cf. \citeauthor{Martino2018}~\cite[in sec. 1.1]{Martino2018}). 

Before constructing the computer, the numerical analysis of deterministic models used approximations that aimed, among other things, at reducing and simplifying the calculations. The ability to perform a significantly larger number of operations in a short time led Ulam (cf. \citeauthor{Met1987:BeginMC}~\cite{Met1987:BeginMC}) to the idea of transforming a deterministic task to an appropriate stochastic task and applying the numerical analysis using the simulation of random quantities of the equivalent model (cf.~\citeauthor{MetUla1949:MCM}~\cite{Met1987:BeginMC},). The problem Ulam was working on was part of a project headed by von Neumann, who accepted the idea, and Metropolis gave it the name of the Monte Carlo Method \cite[in sec. 1.1]{Martino2018}.

The resulting idea allowed us to free oneself from difficult calculations, replacing them with a large number of easier calculations. The problem was the errors that could not be solved by increasing the number of iterations. Work on improving the method consisted in reducing the variance of simulated samples, which led to the development of the Importance Sampling (IS) method in 1956 (cf. \citeauthor{Mar1956:IS}~\cite{Mar1956:IS} ).

Computational algorithms in this direction are a rapidly growing field. Currently, the computing power of computers is already large enough, to provide an opportunity to solve problems that can't be solved analytically. However, some problems and methods still require a lot of power. A lot of research is devoted to finding or improving such methods. The development of this field is currently very important. There are many models describing current problems. Not all offer the possibility of an analytical solution. That is why various studies are appearing to find best methods for accurate results. the Cross-Entropy method helps to realize \textbf{IS} and it may, among other things, be an alternative to current methods in problems in epidemiological models.

\subsection{Sequential Monte Carlo methods.}		
The \textbf{CE} method was proposed as an algorithm for estimation of  probabilities of rare events.  Modifications of cross entropy to minimize variance were used for this by \citeauthor{rubinstein1997optimization} in his seminal paper \cite{rubinstein1997optimization} in \citeyear{rubinstein1997optimization}. This is done by translating the "deterministic" optimization problem in the related "stochastic" optimization and then simulating a rare event. However, to determine this low probability well, the system would need a large sample and a long time to simulate. 

\subsection{Importance sampling.}
Therefore, it was decided to apply the importance sampling (\textbf{IS})\label{IS::important}, which is a technique used for reducing the variance. The system is simulated using other parameter sets (more precisely, a different probability distribution) that helps increase the likelihood of this occurrence. Optimal parameters, however, are very often difficult to obtain. This is where the important advantage of the \textbf{CE} method comes in handy, which is a simple procedure for estimating these optimal parameters using an advanced simulation theory. \citeauthor{sani2007optimal} \cite{sani2007optimal} have used the Cross-Entropy method as the support of the importance sampling to solve the problem of the optimal epidemic intervention of HIV spread. The idea of his paper will be adopted to treat the \textbf{SIRC} models, which has a wide application to modeling the division of population to four groups of society members based on the resistance to infections (v. modeling of bacterial Meningitis transmission dynamics analyzed by \citeauthor{AsoNya2018:Meningitis}~\cite{AsoNya2018:Meningitis}, \citeauthor{Ver2008:SIRC}~\cite{Ver2008:SIRC}).   The result will be compared with an analytical solution of \citeauthor{casagrandi2006sirc}~\cite{casagrandi2006sirc}. 

\subsection{\label{section:cemethod}Cross-Entropy method.} This part will contain information on the methods used in the paper. The Cross-Entropy method developed in \cite{rubinstein1997optimization} is one of the versions of the Monte Carlo method developed for problems requiring the estimation of events with low probabilities. It is an alternative approach to combinatorial and continuous, multi-step, optimization tasks. In the field of rare event simulation, this method is used in conjunction with weighted sampling, a known technique of variance reduction, in which the system is simulated as part of another set of parameters, called reference parameters, to increase the likelihood of a rare event. The advantage of the relative entropy method is that it provides a simple and fast procedure for estimating optimal reference parameters in IS. 

Cross entropy is the term from the information theory (v. \cite{Par1969:Entropy}, \cite[Chapter 2]{CovTho2006:ElementIT}). It is a consequence of measuring the distance of the random variables based on the information provided by their observation. The relative entropy, the close idea to the cross entropy,  was first defined by  Kullback  and  Leibler~\cite{KulLei1951:Information} as a measure of the distance  $\mathcal{D}(p||q)=\textbf{E}_p\log\frac{p(X)}{q(X)}$ (v. \cite[(1.6) on p. 9; Definition on Sec. 2.3  p. 19]{CovTho2006:ElementIT}, \cite{Ing2014:TIandSM}) for two distributions $p$, $q$  and is known also as Kullback-Leibler distance. This idea was studied in detail by Csisz\'ar \cite{Csi1967:Information} and Amari~\cite{Ama1985:Diff}. The application of this distance measure in the Monte Carlo refinements is related to a realization of the importance sampling technique. 

The \textbf{CE} algorithm can be seen as a self-learning code covering the following two iterative phases:
	\begin{enumerate}\itemsep1pt \parskip1pt \parsep1pt
		\item Generating random data samples (trajectories, vectors, etc.) according to a specific random mechanism.
		\item Updating parameters of the random mechanism based on data to obtain a “better” sample in the next iteration.
	\end{enumerate}
Now the main general algorithm behind the Cross-Entropy method will be presented on the examples. 

\subsection{Application to optimization problems.} This section is intended to show how extensive the \textbf{CE} method is. The examples will show the transformation of the optimization problems to the Monte Carlo task of estimation for some expected values. The equivalent \textbf{MC} task uses the \textbf{CE}  method. Other very good examples can be found in \cite{Rub1999:CEMethod}.
    
The first variational problem presented in this context is the specific multiple stopping model. Consider the so-called secretary problem (v. Ferguson~\cite{Fer1989:Who} or the second author's paper~\cite{sza1982:Ath}) for multiple choice, i.e. the issue of selecting the best proposals from a finite set with at most $k$ attempts. There is a set with $N$ objects numbered from $1$ to $N$. By convention, the object with the number 1 is classified as the best and the object with the number $N$ as the worst. Objects arrive one at a time in a random order. We can accept this object or reject it. However, if we reject an object, we cannot return. The task is to find such an optimal stopping rule that the sum of the ranks of all selected objects is the lowest (their value is then the highest). So we strive for a designation
	\[
	\inf_{\tau}\bE X_{\tau}=\bE X_{\tau^{*}}.
	\]
Let $a_{i}$ be the rank of the selected object. The goal is to find a routine for which the value of $\bE (a_{\tau_{1}}+\dots+a_{\tau_{k}})$ for $k\geqslant2$, is the smallest. More details are moved to~\ref{SumOfRank}.

Next, the problem was transformed to the minimization of the mean of sum of ranks. The details are presented in the section~\ref{MinMeanSumRank}.  This problem and its solution by \textbf{CE} was described by \citeauthor{Pol2010:CEBCP}~\cite{Pol2010:CEBCP}. Its correctness was checked and programmed in a different environment by \citeauthor{Sta2019} in \cite{Sta2019}. 
	
The second example is presented in details in Section~\ref{VehRP}. It is a formulation of the vehicle routing problem (v. \cite{Dro2002:VRSD}, \cite{CheHom2005:vehicle}). This is an example of stochastic optimization where the Cross-Entropy is used.

\subsection{Goal and organization of this paper. }
In the following parts of the paper, we will focus on models of the dynamics of the spread of infection over time when the population consists of individuals susceptible to infection, sick, immune to vaccination or past infection, and partially susceptible i.e. the population is assigned to four compartments of individuals. Actions taken and their impact on the dynamics of the population are important. It is precisely the analysis of the impact of the preventive diagnosis and treatment that is particularly interesting in the model. How the mathematical model covers such actions is presented in Section~\ref{section:sircmodel}. The model created in this way is then adapted for Monte Carlo analysis in~\ref{section:controlfunction} and the results obtained on this basis are found in~\ref{section:applicationcemethod}. 
The analysis of computational experiments concludes this discussion in~\ref{chapter:results}.


\section{\label{chapter:sircwithcemethod} Optimization of control for \textbf{SIRC} model.}
Let us focus the attention on the introduction of \textbf{SIRC} model and presentation of the logic behind using the \textbf{CE} method to determine functions that optimize the spread of epidemics. The ideology behind creation of the \textbf{SIRC} model and its interpretation will be presented in order to match the right parameters for calculating the cost functions. 
	
	The \textbf{CE} method here is used to solve the variational deterministic problem. Solving such problems has long been undertaken with the help of numerical methods with a positive result. In the position \cite{ekeland1999convex} proposals of the route for variational problems for optimization problems in time are presented. The main focus was on the non-convex problems, which often occur with optimal control problems. In the book by Glowinski~\cite{glowinski2008lectures} a review of the methods for solving variational problems was made. And then \cite{mumford1989optimal} presented a complicated method of approximating functions for the problem optimization.

\subsection{\label{section:sircmodel}\textbf{SIRC} model.} The subject of the work is to solve the problem optimizing the spread of the disease, which can be modeled with the \textbf{SIRC} model. This model was proposed by \citeauthor{casagrandi2006sirc}
	\cite{casagrandi2006sirc}. Its creation was intended to create a better model that would describe the course of influenza type A. Contrary to appearances, this is an important topic because the disease, although it is widely regarded as weak, is a huge problem for healthcare. In the US alone, the cost associated with the influenza epidemic in the season is estimated to exceed USD 10 billion (v.~\cite{klimov1999surveillance}) and the number of deaths is over 21,000 per season (v.~\cite{simonsen1997impact}).

	Various articles have previously proposed many different mathematical models to describe a pandemic of influenza type A. An overview of such models was made by Earn~\cite{earn2002ecology}. In general, it came down to combining SIR models with the use of cross-resistance parameters (cf. details in the papers by Andreasen et al.~\cite{AndLinLev1997:Dynamics} and Lin et al.~\cite{lin1999dynamics}). The authors of the article write that the main disadvantage of this approach is the difficulty of the analysis and calculations with a large number of strains. The results of their research showed that the classic SIR model cannot be used to model and study influenza epidemics. Instead, they proposed the extension of the SIR model to a new class C that will simulate a state between a susceptible state and a fully protected state. This extension aims to cover situations where vulnerable carriers are exposed to similar stresses as they had before. As a result, their immune system is stronger when fighting this disease. 

\begin{figure}[H]
		\centering
		\includegraphics[width=0.6\linewidth]{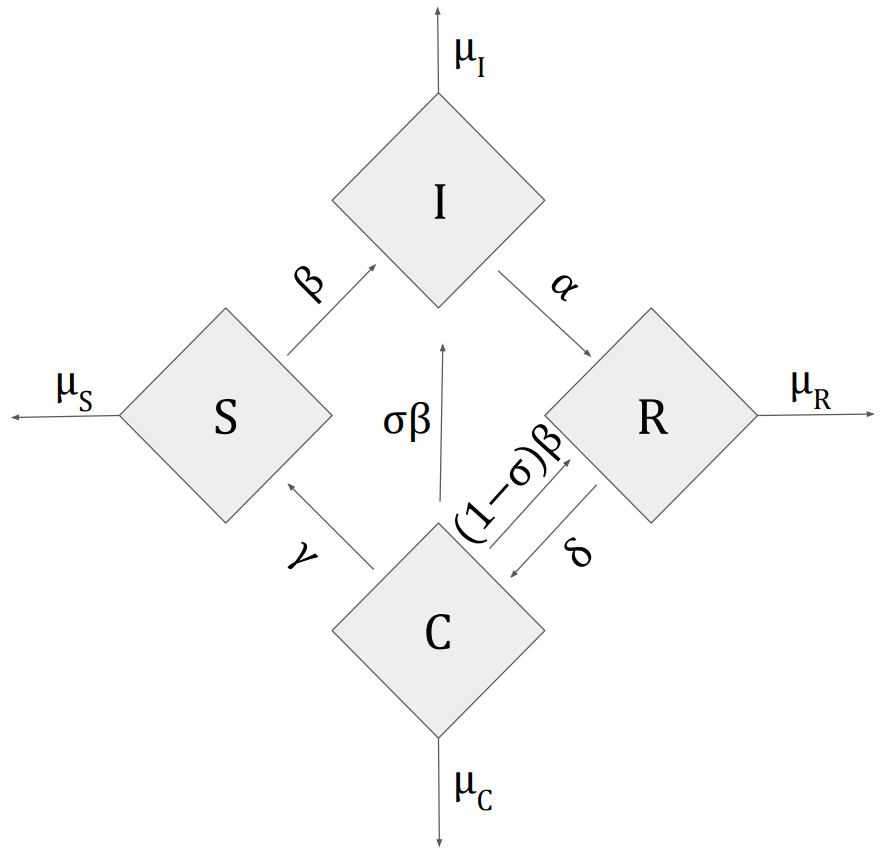}
		\caption{\label{fig:sirc1}Schematic representation of the \textbf{SIRC} model.}		
	\end{figure}

The \textbf{SIRC} model divides the community into 4 groups:
	\begin{tableitems}\itemsep2pt \parskip3pt \parsep2pt
		\item S - persons susceptible to infection, who have not previously had contact with this disease and have no immune defense against this strain
		\item I - people infected with the current disease
		\item R - people who have had this disease and are completely immune to this strain
		\item C - people partially resistant to the current strain (e.g. vaccinated or those who have had a different strain)
	\end{tableitems}
	Figure \ref{fig:sirc1} contains a general scheme of this model. The four rhombus represent the four compartments of individuals, the movement between the compartments is indicated by the continuous arrows. The main advantage of this model over other SIR users is the fact that after recovering from a given strain, in addition to being completely resistant to this strain, they are also partly immune to a new virus that will appear later. This allows you to model the resistance and response of people in the group to different types of disease.
     Parameters	$\alpha, \delta$ i $\gamma$ can be interpreted as the reciprocal of the average time spent by a person in order of ranges $I, R, C$. Parameters $\mu$ with indices $S,I,R,C $ represent the natural mortality in each group, respectively. In some versions of the model, additional mortality rates are considered for the group of infected people. Here we assume that this factor is not affected by the disease and we will denote it as $\mu$ in later formulas ($\mu_{S} = \mu_{I} = \mu_{R} =\mu_{C} = \mu$). The next one $\sigma$ is the likelihood of reinfection of a person who has cross-resistance while the parameter $\beta$ describes the contact indicator.
	
	The \textbf{SIRC} model is represented as a set of four ordinary differential equations. Let $S(t)$, $I(t)$ be the number of people in adequate compartments. We have (v. \cite{casagrandi2006sirc})
	\begin{equation} \label{eq:ode_1}
		\begin{aligned}	 
		\frac{dS}{dt} &= \mu(1-S(t))-\beta S(t)I(t)+ \gamma C(t), \\
		\frac{dI}{dt} &= \beta S(t)I(t) + \sigma \beta C(t)I(t)-(\mu+\alpha)I(t), \\
		\frac{dR}{dt} &= (1- \sigma) \beta C(t)I(t) + \alpha I(t) - (\mu+\delta)R(t),\\
		\frac{dC}{dt} &= \delta R(t) - \beta C(t)I(t) - (\mu+\gamma)C(t),
		\end{aligned}
	\end{equation}
	with initial conditions
	\begin{equation} \label{eq:ode_2}
		\begin{split}	 
		S(0)=S_{0}, \, I(0)=I_{0}, \, R(0)=R_{0}, \, C(0)=C_{0} 
		\end{split}
	\end{equation}
	
	\subsection{\label{section:controlfunction}Derivation of optimization functions}		
An approach to optimize this model was proposed by~\citeauthor{iacoviello2013optimal}~\cite{iacoviello2013optimal}. In this article one can find the suggestion concerning performing the calculations as outlined by \citeauthor{zaman2008stability} in \cite{zaman2008stability}. Two functions were proposed in the approach as the parameter or controls. One relates to the \emph{susceptible people} and the other describes the number of \emph{sick people}. The method described by \citeauthor{kamien2012dynamic}~\cite{kamien2012dynamic} was used to determine the optimal solution. In order to apply it, the set of equations~\eqref{eq:ode_1} should be updated accordingly. After these modifications, it has the following form: 
	\begin{equation} \label{eq:ode_3}
	\begin{aligned}	 
		\frac{dS}{dt} &= \mu(1-S)-\beta SI+ \gamma C - g(S(t),u(t)), \\
		\frac{dI}{dt} &= \beta S(t)I(t) + \sigma \beta C(t)I(t)-(\mu+\alpha)I(t) - h(I(t),v(t)), \\
		\frac{dR}{dt} &= (1- \sigma) \beta C(t)I(t) + \alpha I(t) - (\mu+\delta)R(t) + g(S(t),u(t)) + h(I(t),v(t)),\\
		\frac{dC}{dt} &= \delta R(t) - \beta C(t)I(t) - (\mu+\gamma)C(t),
		\end{aligned}
	\end{equation}
	where
	\begin{equation} \label{eq:ode_4}
		\begin{aligned}	 
	g(S(t),u(t))&=\rho_{1}S(t)u(t),\\
	h(I(t),v(t))&=\rho_{2}I(t)v(t).
		\end{aligned}
	\end{equation}
	$u(t)$ represent the percentage of those susceptible who have been taken care of thanks to using control on population and  $v(t)$ represent the percentage of infected people with the same description, respectively. $\rho_{1}$ and $\rho_{2}$ are weights that optimize the proportion of given control options. Both functions $g(S(t),u(t))$ $h(I(t),v(t))$ and in order they can be interpreted as actions performed on people susceptible to the disease and infected. In addition, two new conditions related to optimization functions are added to the initial conditions.
	\begin{equation} \label{eq:ode_5}
	\begin{split}	 
	u_{min} \leqslant u(t) \leqslant u_{max}, \, v_{min} \leqslant v(t) \leqslant v_{max}
	\end{split}
	\end{equation}
	They describe the limits on what part of the population care can be given to. The smallest values $u_{min}$ and $v_{min}$ is zero.
	The existence of the solution is shown by Iacoviello et al.~\cite{iacoviello2013optimal} and they are functions that minimize the following cost index:
	\begin{equation} \label{eq:ode_6}
	J(u,v) = \int_{t_{1}}^{t_{2}}  \varphi(S(t),I(t),u(t),v(t))  dt 
	= \int_{t_{1}}^{t_{2}} \left[ \alpha_{1}S(t) + \alpha_{2}I(t) + \frac{1}{2}\tau_{1}u^{2}(t) + \frac{1}{2}\tau_{2}v^{2}(t) \right] dt
	\end{equation}
	$\alpha_{1},\alpha_{2}$ are used to maintain the balance in the susceptible and infected group. $\tau_{1},\tau_{2}$ is interpreted as weighting in the cost index. $t_ {1}, t_ {2}$ determine the time interval. The square next to the functions $u(t)$ and $v(t)$ indicates increasing intensity of functions \cite{zaman2008stability}. Thus the function $J(u,v)$ reflects the human value of susceptible and infected persons, taking into account the growing value of funds used over a specified period of time.

    The objective function can be changed as long as there is a minimum. The \textbf{CE} method does not impose any restrictions here. However, to compare the results with \cite{iacoviello2013optimal}, the function proposed by them is used also here. There, the objective function must be square in order for the quadratic programming techniques to be used. Usually quadratic function is good in mathematics, but not in practice. Changing the objective function and impact of it can be considered further.
    
	\subsection{\label{section:applicationcemethod}Optimization of the epistemological model by  \textbf{CE} method.} 
	\citeauthor{sani2007optimal} \cite{sani2007optimal} present the way in which the Cross-Entropy method can be used to solve the optimization of the epistemological model consisting of ordinary differential equations. The main task was to minimize the objective function $J(u)$ depending on one optimization function $u(t)$ over a certain set $U$ consisting of continuous functions $u$. The minimum can be saved as:
	\begin{equation}\label{eq:ce_1}
	\gamma^{*} =J(u^{*})=\min _{u\in U}J(u)
	\end{equation}
	Parameterizing the minimum problem looks as follows:
	\begin{equation}\label{eq:ce_2}
	\min _{c\in C}J(u_{c})
	\end{equation}
	where $u_{c}$ is a function from $U$, which is parameterized by a certain control vector $c\in \mathbb{R}^{m}$. Collection $C$ is a set of vectors $c$. It should be noted that the selected set $C$ should be big enough to get the enough precise solution. Then in such a set there are such $c^{*}$-- which will be the optimal control vector and  $\gamma^{*}$--optimal value, respectively.
	 
    One way to parameterize a problem is to divide the time interval $[t_{1},t_{2}]$ at small intervals and using these intervals together with control vector points $c$ to define a function $u_{c}$. This function can be created by interpolating between points. Such interpolation can be done using e.g. the finite element method, finite difference method, finite volume method or cubic B-Spline . \citeauthor{fang2002finite} in \cite{fang2002finite} and \citeauthor{caglar2006} in \cite{caglar2006} show that all of this methods can be used to the two-point boundary value problems.  Here the \textbf{FEM}\label{FEM} method was used (v. \cite{PapGio2018:SFEM}). The results for the other methods have not been checked. The choice of the approximation method is a difficult and interesting issue, but this is not what is consider here. With this function $u_{c}$, it takes form 
	\begin{equation}\label{eq:ce_3}
	u_{c}(t) = \sum_{t=0}^{n} c_{i}k_{i}(t)
	\end{equation}
	where
	\begin{equation}\label{eq:ce_4}
	  k_{i}(t) =
	\begin{cases}
	\frac{t-t_{i-1}}{t_{i}-t_{i-1}} & \text{for $t \in [t_{i-1},t_{i}]$}\\
	\frac{t_{i+1}-t}{t_{i+1}-t_{i}} & \text{for $t \in [t_{i},t_{i+1}]$}\\
	0 & \text{otherwise.}
	\end{cases}
	\end{equation}
	
	Now when $u_{c}$ is represented as a parameterized function, the value can be counted$J(u_{c})$ by solving the system of ordinary differential equations \eqref{eq:ode_3} e.g. by the Runge-Kutta method. The idea of using the \textbf{CE} method in this problem is to use a multi-level algorithm with generated checkpoints $c_{i}$, here, from a normal distribution. The \textbf{CE} method does not impose the distribution. Usually it is chosen from the family of density which is expected. Here, it's expected that there will be one high peak at the peak of the epidemic. Hence the normal distribution. Then, update the distribution parameters using the \textbf{CE} procedure using the target indicator as a condition for selecting the best samples. The calculations are carried out until the  empirical variance of optimal control function is smaller than the given $\epsilon$. This means that the function values are close to the optimal expected value  The idea of applying the \textbf{CE} method to the \textbf{SIRC} model is similar. The only difficulty is the fact that when optimizing the \textbf{SIRC} model there are two functions instead of one as in the case of the algorithm described here.
	\noindent In the form of an algorithm, it looks as follows: \\[2ex]
	Algorithm ( \emph{Modification of the \citeauthor{sani2007optimal}'s algorthm (v. ~\cite{sani2007optimal}) to two optimal functions case}):
	\begin{enumerate}\itemsep3pt \parskip3pt \parsep2pt
		\item Initialize: Choose $\mu^{(u)}_{0} = \{\mu_{i0}, i=0,\dots,n\}$, $\sigma^{(u)}_{0}=\{ \sigma_{i0},i=0,\dots,n \}$, $\mu^{(v)}_{0} = \{\mu_{i0}, i=0,\dots,n\}$ and $\sigma^{(v)}_{0}=\{ \sigma_{i0},i=0,\dots,n \}$, set $k=1$
		
		\item Draw: Generate a random samples $C_{1},\dots,C_{N} \sim N(\mu^{(u)}_{k-1},{\sigma^{(u)}_{k-1}}^{2})$ and $D_{1},\dots,D_{N} \sim N(\mu^{(v)}_{k-1},{\sigma^{(v)}_{k-1}}^{2})$ with $C_{m}=\{ C_{mi}, i=0,\dots,n \}$ and $D_{l}=\{ D_{li}, i=0,\dots,n \}$.\\[2ex]
		
		\item\label{AlgEvalJ1} Evaluate: For each control vector $C_{m}$ and $D_{m}$ evaluate the objective function $J(u_{C_{m}},u_{D_{l}})$ by solving the ODE system
		
		\item\label{AlgEvalJ2} Find the $p$ best performing samples, based on the values $\{J(u_{C_{m},D_{l}})\}$. Let $I$ be the corresponding set of indices.
		
		\item Update: for all $i=0,\dots,n$ let
		\begin{equation*}
		\begin{aligned}
		\hat{\mu_{ki}}^{(u)} = \frac{1}{p} \sum_{m\in I}^{}C_{mi};\qquad 
		(\hat{\sigma_{ki}}^{(u)})^{2} = \frac{1}{p} \sum_{m\in I}^{}(C_{mi}-\mu^{(u)}_{ki}), \\
		\hat{\mu_{ki}}^{(v)} = \frac{1}{p} \sum_{l\in I}^{}D_{li}; \qquad
		(\hat{\sigma_{ki}}^{(v)})^{2} = \frac{1}{p} \sum_{l\in I}^{}(D_{li}-\mu^{(v)}_{ki}).
		\end{aligned}
		\end{equation*}
		
		\item Smooth: For a fixed smoothing parameter $ 0 < \alpha \leqslant 1$ let
		\begin{equation*}
		\begin{aligned}
		\hat{\mu_{k}}^{(u)} &= \alpha\mu_{k}^{(u)}+(1-\alpha)\hat{\mu_{k-1}}^{(u)};\qquad 
		\hat{\sigma_{k}}^{(u)} &= \alpha\sigma_{k}^{(u)}+(1-\alpha)\hat{\sigma_{k-1}}^{(u)} \\
		\hat{\mu_{k}}^{(v)} &= \alpha\mu_{k}^{(v)}+(1-\alpha)\hat{\mu_{k-1}}^{(v)}; \qquad
		\hat{\sigma_{k}}^{(v)} &= \alpha\sigma_{k}^{(v)}+(1-\alpha)\hat{\sigma_{k-1}}^{(v)}
		\end{aligned}
		\end{equation*}		
		
		\item Repeat 2–6 until $max_{i}\sigma_{ki}^{(u)}< \epsilon$ and $max_{i}\sigma_{ki}^{(v)}< \epsilon$ with $\epsilon = 10^{-5} $ . Let $L$ be the final iteration number. Return $\mu_{L}$ as an estimate of the optimal control parameter $c^{*}$
	\end{enumerate}
    The value of the  criterion $ J(\cdot, \cdot) $ is calculated in steps \ref{AlgEvalJ1} and \ref{AlgEvalJ2}. The selection of the function $\varphi$ allows to fit the cost function to the modeled case. Here, the criterion given by \eqref{eq:ode_6} is adopted.
	

	\section{\label{chapter:results}Description of the numerical results.}	
	This section will present the results obtained using the Cross-Entropy method. For the correct comparison the same parameters for the \textbf{SIRC} model were used. These values are:
	\begin{equation*}
	\begin{aligned}
	\mu&=\frac{1}{75}[\text{year}^{-1}], \, &\gamma&=\frac{1}{2}[\text{year}^{-1}], \, &\alpha&=\frac{365}{5}[\text{year}^{-1}], \\
	\delta&=1[\text{year}^{-1}], \, 
	&\beta &\approx 146 [\text{year}^{-1}], \, &\sigma&\approx 0.078.
	\end{aligned}
	\end{equation*} 
	For optimization functions $u(t)$ and $v(t)$ the following restrictions have been applied:
	\begin{equation*} 	 
	0 \leqslant u(t) \leqslant 0.9, \, 0 \leqslant v(t) \leqslant 0.9
	\end{equation*}
	Restrictions were proposed by \citeauthor{lenhart2007optimal}~\cite{lenhart2007optimal}). Their value was explained by the fact that the whole group cannot be controlled. In addition, the weights used for functions have values $\rho_{1}=2$ and $\rho_{2}=2$ \\
	Now, all that's left is to propose variable values for the objective function:
	\begin{equation*} 	 
	\alpha_{1}=10^{-3}, \, \alpha_{2}=0.997, \, \tau_{1}=10^{-3}, \, \tau_{2}=10^{-3}.
	\end{equation*}
	
	The parameter values are adjusted using historical data like parameter mortality rate $\mu$, which is counted as average lifetime of a host or parameters $\alpha, \delta, \gamma$ mean the inverted time of belonging to the group $I, R$ and $C$. A detailed description of determining parameters along with their limitations is presented in the article by \citeauthor{casagrandi2006sirc}~\cite{casagrandi2006sirc}.
	
    \subsection{A remark about adjusting the parameters of the control determination procedure.} Here these parameters are given and describes the influenza A epidemics well, because of long historical data sets. However, this may not be the case, and then these parameters are subject to adjustment in the control determination procedure. An overview of applications and accuracy of calibration methods, which can be used for it is presented in the article by \citeauthor{hazelbag2020calibration}~\cite{hazelbag2020calibration}, which provides an overview of the model calibration methods. All parameters in the model are constant and independent of time, so methods which try to optimizes a goodness-of-fit (\textbf{GOF}\label{GOF}) can be use to this example. \textbf{GOF}\label{GOF} is a measure that assesses consistency between model the output and goals. As a result, it gives the best combination of parameters. Examples of such methods are Grid search, Nelder-Mead method, the Iterative, descent-guided optimisation algorithm, Sampling from a tolerable range. After finding the appropriate parameters, other algorithms like the profile likelihood method or Fisher information can be used to calculate the confidence intervals for these coefficients. If the epidemic described by the model consists of transition probabilities that cannot be estimated from currently available data, calibrations can be performed to many
	end points. Then \textbf{GOF} is measured as the mean percentage deviation of the values obtained at the endpoints (v. \cite{taylor2010methods}).
	
	\subsection{Proposed optimization methods in the model analysis.}
\begin{figure}[tbh!]
		\centering
		\includegraphics[width=0.75\linewidth]{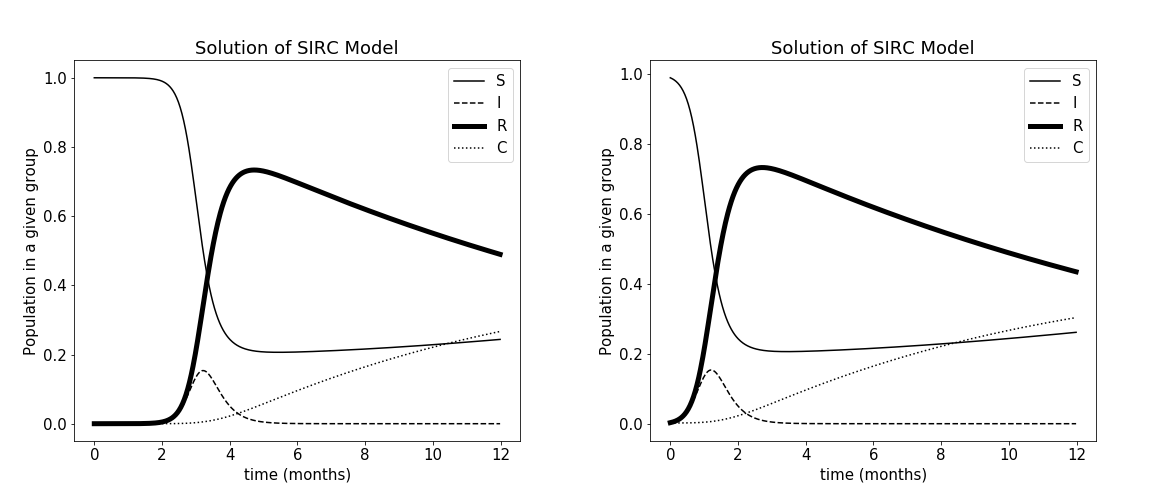}
		\caption{\label{fig:sircsolutionwithout}Solutions for \textbf{SIRC} models without optimal control (in order the epidemic started and developed)}
	\end{figure}

	The results will be presented for two moments in the model: at the beginning and in the middle of the epidemic. This is initialized with other initial parameters. For the beginning of the epidemic, they are as follows:
	\begin{equation*} 	 
	S_{0} = 1-I_{0}, \, I_{0}=10^{-6}, \, R_{0}=0, \, C_{0}=0.
	\end{equation*}
	And for the widespread epidemic:
	\begin{equation*} 	 
	S_{0} = 0.99, \, I_{0}=5*10^{-3}, \, R_{0}=3*10^{-3}, \, C_{0}=2*10^{-3}.
	\end{equation*}
	The values in the model have been normalized for the entire population, and they add up to 1. The results for both models without the use of controls are shown in \ref{fig:sircsolutionwithout}. They were obtained using the Runge-Kutty method. The cost index in this case was in order $0.00799$ and $0.00789$ for the model of the beginning and the development of an epidemic.

In Sections~\ref{section:version1} and \ref{section:version2} two ways of calculating optimal functions $u(t)$ and $v(t)$ for the \textbf{SIRC} model will be considered. In the first version, the functions will be calculated separately. In the next both $u(t)$ and $v(t)$ will be counted at the same time. All versions will be presented at two moments of the epidemic: when it began and when it has already spread. The results will be compared with those obtained in the article by Iacoviello et al.~\cite{iacoviello2013optimal}, where the problem was solved using the sequential quadratic programming method using a tool from Matlab. Unfortunately, the article does not specify the cost index value for solving the obtained sequential quadratic programming method. Therefore, in the paper, there were attempts to recreate the form of the function $u(t)$ and $v(t)$ and the cost index obtained for them ($0.003308$ for the first situation and $0.006489$ for the second one) and the model solution.

\begin{figure}[tbh!]
		\centering
		\includegraphics[width=.75\linewidth]{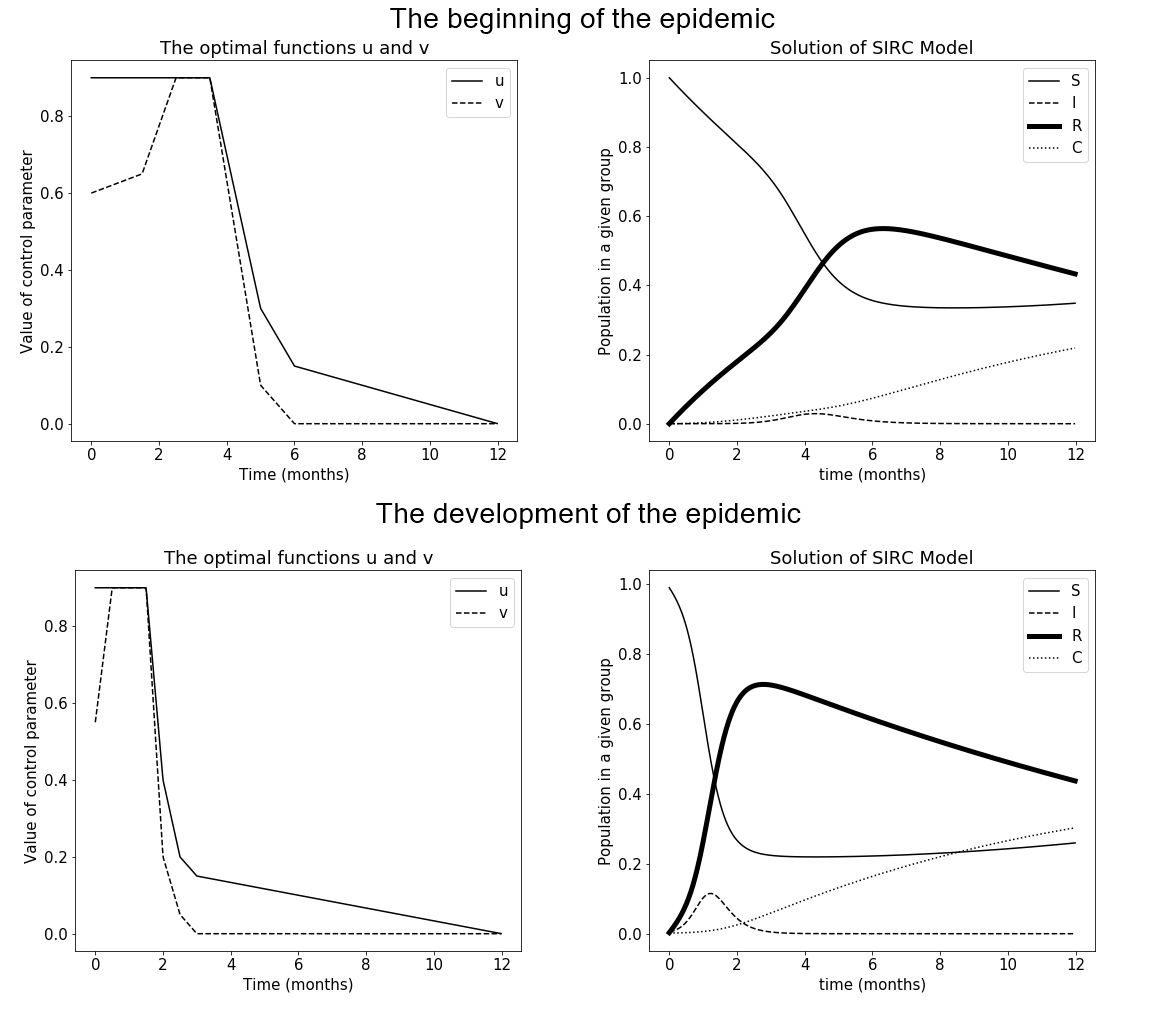}
		\caption{Values of control functions and solutions obtained for them. Solved by the sequential quadratic programming method}
		\label{fig:ceversion0}
	\end{figure}

	As can be seen in Figure \ref{fig:ceversion0} in the first situation, the control functions helped reduce peak infection, which was previously without control (Figure \ref{fig:sircsolutionwithout}) occurred between 2 and 4 months. Previously, the largest peak was around $20\%$ occurred between 2 and 4 months. Previously, the largest peak was around $10\%$. In the case of the second situation, the peak also decreased, but not very clearly (it is also visible in the value of the cost index, which for the first situation is $0.003308$, and for the other one $0.006489$). It can be seen that timely control is also very important to reduce the harmful effects of an epidemic. For more conclusions on how to apply control properly in this model, see \cite{iacoviello2013optimal}.	

	\subsection{\label{section:version1}\textbf{CE} method version 1.}\unskip In the first method, the algorithm shown at the end of the section \ref{section:applicationcemethod} has been applied twice. Once for the calculation of the optimal function $u(t)$ with a fixed $v(t)$, and then $v(t)$ with a fixed $u(t)$. In the next step, the results were combined and the result for the \textbf{SIRC} model was calculated using these two functions. This approach is possible due to the lack of a combined restriction on $u(t)$ and $v(t)$. Features are not dependent on each other with resources.

\begin{figure}[tbh!]
		\centering
		\includegraphics[width=.7\linewidth]{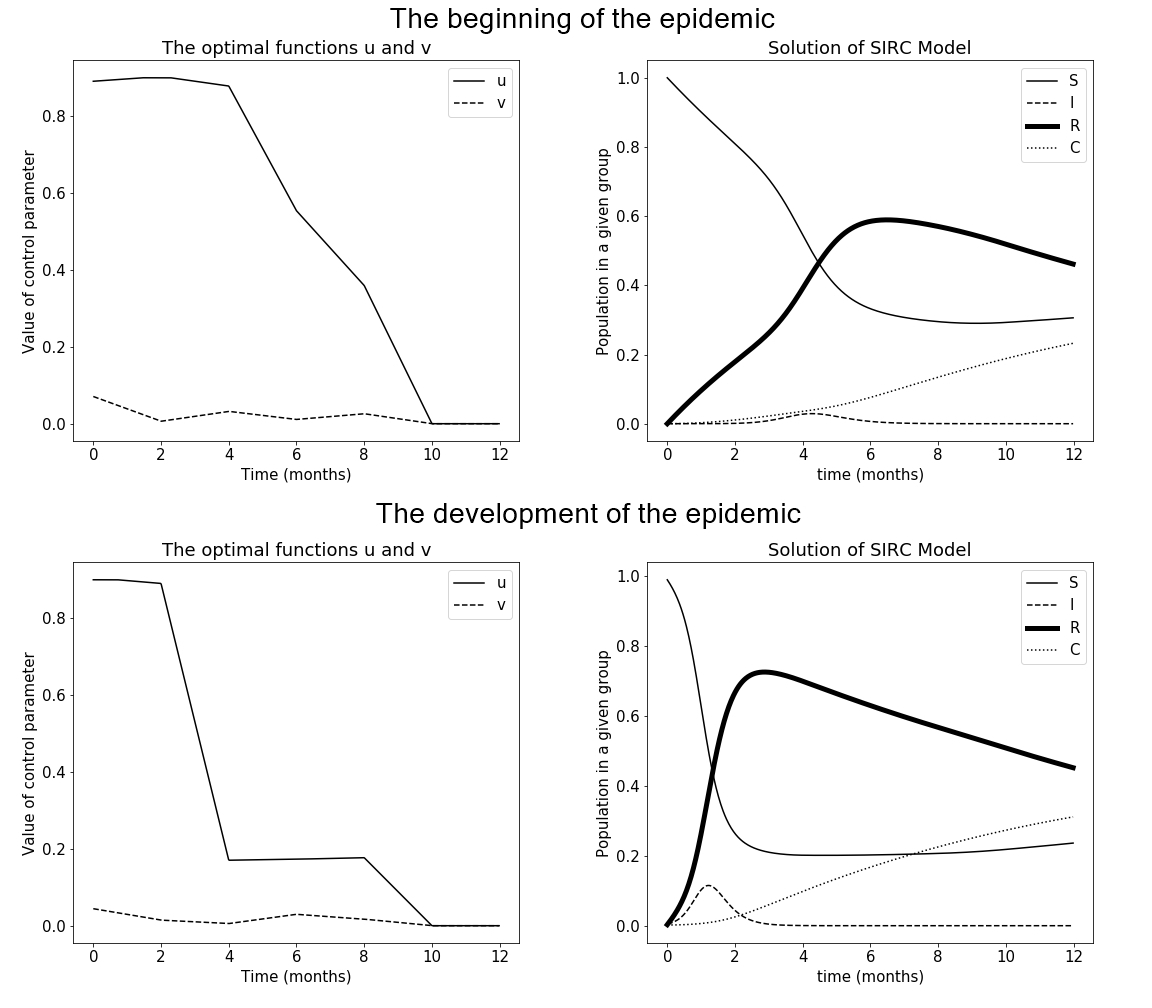}
		\caption{Values of control functions and solutions obtained for them. Solved by the cross-entropy method - version 1}
		\label{fig:ceversion1}
	\end{figure}

	The cost index was $0.003086$ for the start of epidemic and $0.006443$ a developed epidemic. The figure~\ref{fig:ceversion1} show the graphs of functions $u(t)$ and $v(t)$  along with a graph of solutions for both moments during the epidemic.

	For functions $u(t)$ the results are similar, only the decrease here is more linear. However, the solution of the \textbf{SIRC} model gives the same results. Function $v(t)$ is much more interesting. The \textbf{CE} method showed that it should have very low values, even close to $0$ throughout the entire period in both cases. However, this does not affect the cost index (it is even smaller than in the previous case) and the \textbf{SIRC} model. More on this subject will be included in the summary. 
	
	\subsection{\label{section:version2}\textbf{CE} method version 2.}\unskip 	In the second case, it was necessary to change slightly the algorithm from Section~ \ref{section:applicationcemethod} to be able to use it for two control functions. In point 1 of the algorithm a variable is generated $C_m\sim N(\mu_{k},\sigma_{k}^{2})$, which serves to designate $u_{C_{m}}$, where $\mu_{k},\sigma_{k}$ are requested by the user. Here 2 $\mu$ will need to be initialized $\mu_{k}^{(u)},\mu_{k}^{(v)}$ and 2 $\sigma$: $\sigma_{k}^{(u)},\sigma_{k}^{(v)}$. Then designate $C_{1}^{(u)}, \dots ,C_{N}^{(u)} \sim N(\mu_{k}^{(u)},\sigma_{k}^{(u)2})$ and similarly $C_{1}^{(v)}, \dots ,C_{N}^{(v)} \sim N(\mu_{k}^{(v)},\sigma_{k}^{(v)2})$. Now $J(u,v)$ can be count without problems. Further the algorithm goes the same way, and in point 6 all 4 parameters $\mu_{k}^{(u)},\mu_{k}^{(v)},\sigma_{k}^{(u)},\sigma_{k}^{(v)}$ are updated. The graphs obtained for these functions are below.
	\begin{figure}[tbh!]
		\centering
		\includegraphics[width=.7\linewidth]{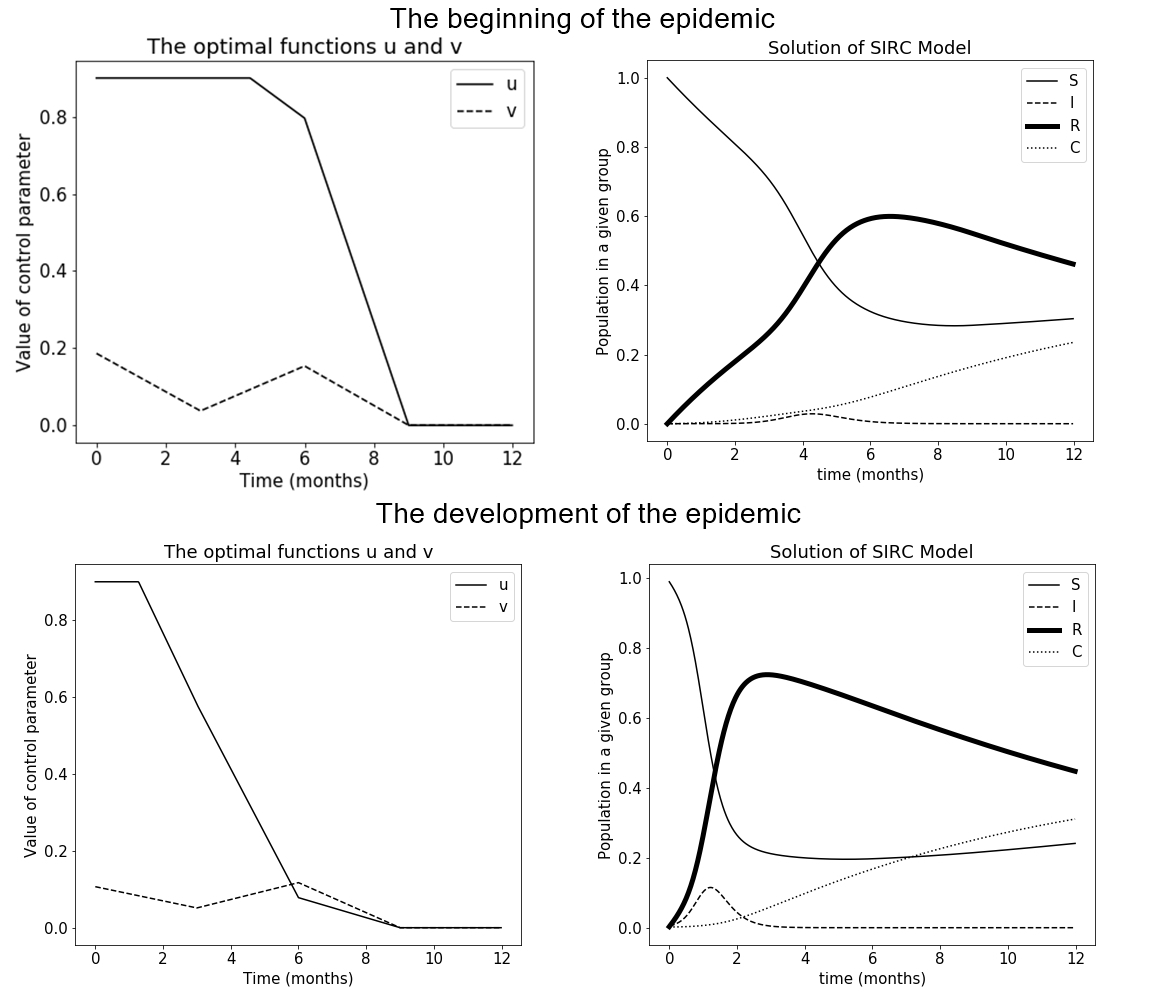}
		\caption{Values of control functions and solutions obtained for them. Solved by the cross-entropy method - version 2}
		\label{fig:ceversion2}
	\end{figure}
	
	In this case, we have very similar results as in the previous section. This can be seen in both the charts and the cost index values ($0.003044$ and $0.006451$). The only problem here is the need for a large sample to bring the results closer. 

\subsection{\label{section:results}Comparison of results.} 	
\begin{figure}[H]
		\centering
		\includegraphics[height=0.25\textheight,width=.95\linewidth]{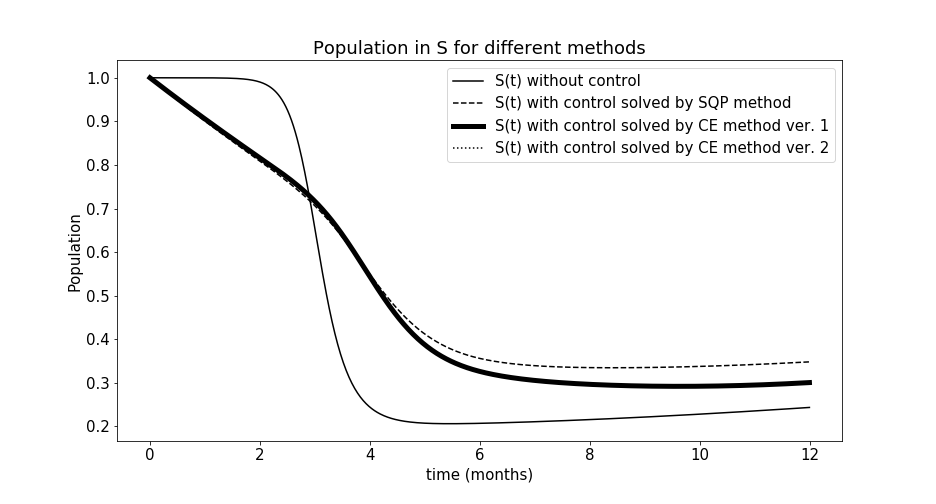}
		\caption{Comparison of results obtained for group S among various methods}
		\label{fig:scomparision}
	\end{figure}

\begin{figure}[H]
		\centering
		\includegraphics[height=0.25\textheight,width=.85\linewidth]{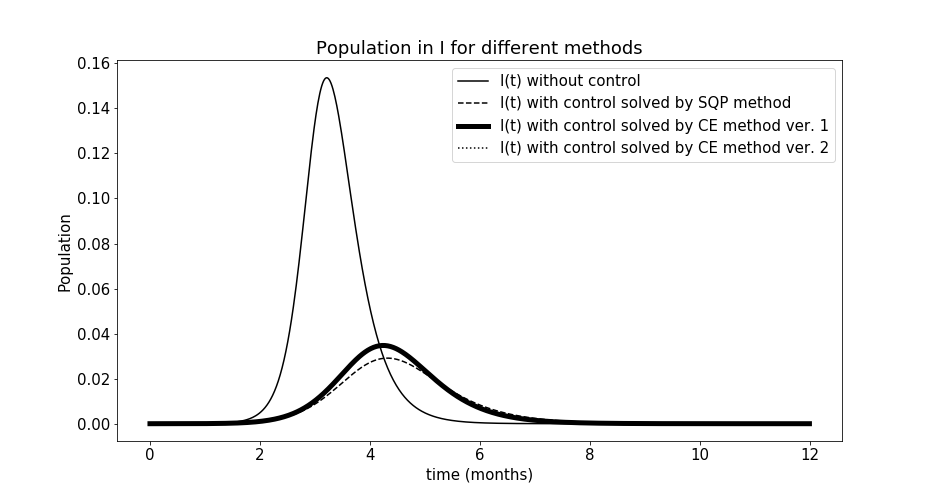}
		\caption{Comparison of results obtained for group I among various methods}
		\label{fig:icomparision}
	\end{figure}

\begin{figure}[H]
		\centering
		\includegraphics[height=0.27\textheight,width=.85\linewidth]{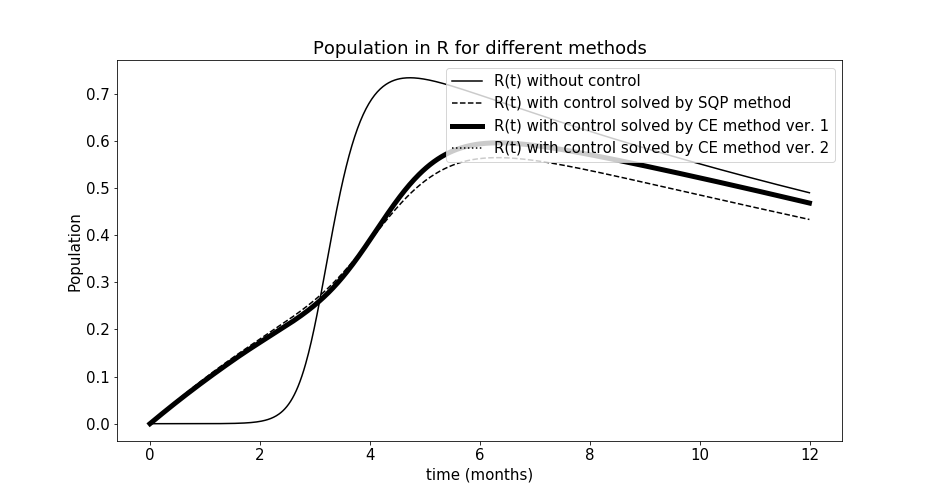}
		\caption{Comparison of results obtained for group R among various methods}
		\label{fig:rcomparision}
	\end{figure}

\begin{figure}[tbh!]
		\centering
		\includegraphics[height=0.27\textheight,width=.85\linewidth]{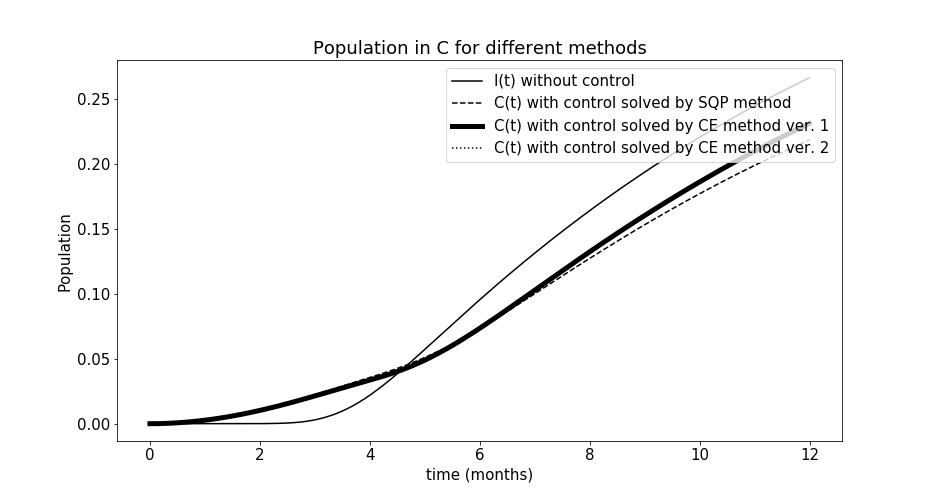}
		\caption{Comparison of results obtained for group C among various methods}
		\label{fig:ccomparision}
	\end{figure}

This section gathers the results from each version and compares them with each other. At first, the behavior of each group of the \textbf{SIRC} model was looked at and compared with the solution of the uncontrolled model. This is shown in Figures \ref{fig:scomparision}, \ref{fig:icomparision}, \ref{fig:rcomparision}, \ref{fig:ccomparision}. The results are presented only for the first situation, where the epidemic begins to spread, because there you can see the results of using control much better.
	
The most interesting is Figure~\ref{fig:icomparision} comparing infected people. You can see here very well how functions $u(t)$ and $v(t)$ help in the fight against the epidemic. Let two versions of the \textbf{CE} method overlap. The functions obtained by the sequential quadratic programming method give slightly better results. However, they have a worse cost index (the comparison is presented in Table \ref{tab:summary}). This is probably due to the high value of the function $v(t)$, which increases the cost index. The question remains why there is such a discrepancy between the obtained function values $v(t)$ with both methods. In this article \cite{iacoviello2013optimal} weight and impact of functions $u(t)$ and $v(t)$ were compared. The chart in the article compares the results obtained for the application of both functions, one at a time and none. It turned out that the main influence on the number of infected has control over susceptible persons, i.e. the function $u(t)$. Choosing the optimal values for $v(t)$ is not very important for the solution of the model. Therefore, such discrepancies are possible when choosing optimal values.	

	The table~\ref{tab:summary} contains a comparison of the cost index values and the duration of the algorithms for each method. The first value is when the epidemic started, and the second is when it has already developed. The article does not specify the algorithm time of the method used.
	\begin{table}[tbh!]
		\centering
		\begin{tabular}{|l|l|l|}
			\hline
			 Method & Cost index & Time \\\hline
			 \multirow{2}{*}{Results without control functions} & \multicolumn{1}{l|}{0.00799} & \multirow{2}{*}{-} \\\cline{2-2}
			 & \multicolumn{1}{l|}{0.00789} & \\\hline			 
			\multirow{2}{*}{Sequential quadratic programming method} & \multicolumn{1}{l|}{0.003308} & \multirow{2}{*}{-} \\\cline{2-2}
			& \multicolumn{1}{l|}{0.006489} & \\\hline				
			\multirow{2}{*}{Cross-entropy method version 1} & \multicolumn{1}{l|}{0.003086} & \multirow{2}{*}{30s} \\\cline{2-2}
			& \multicolumn{1}{l|}{0.006443} & \\\hline				
			\multirow{2}{*}{Cross-entropy method version 2} & \multicolumn{1}{l|}{0.003044} & \multirow{2}{*}{4m54s} \\\cline{2-2}
			& \multicolumn{1}{l|}{0.006451} & \\\hline	
		\end{tabular}
		\caption{\label{tab:summary}Summary of cost indexes and times for all situations}		
	\end{table}
	It can be seen that for each version the cost index values came out very similar. The differences can be seen in time. This is because in the first version of the \textbf{CE} method a smaller sample was simulated because both functions were considered separately. Fortunately, the functions were not closely related and it was possible to consider them separately.

\section{\label{chapter:summary}Summary.}\unskip 
The study examined the possibility of using the \textbf{CE} method to determine the optimal control in the SIRC model. Similar models are being developed to model the propagation of malicious software in computer networks (v. \citeauthor{TayGubZhu2017:OICESHM}~\cite{TayGubZhu2017:OICESHM}). However, the number of works with this approach is still relatively small, although its universality should encourage checking its effectiveness. There are two control functions in the model, which when properly selected optimize the cost function. The \textbf{CE} method was used in two versions: considering the two functions separately and together. The results were then compared with those obtained by \citeauthor{iacoviello2013optimal}~\cite{iacoviello2013optimal} using routing of Matlab\footnote{Matlab is a paid environment with a number of functions. Unlike other environments and programming languages, additional features are created by developers in a closed environment and can only be obtained through purchase.}. The Cross-Entropy method had little problems when considering two functions at the same time. Due to the large number of possibilities, it was necessary to simulate a large sample, which significantly extended the algorithm's time. The second way came out much better. However, it only worked because the functions were not dependent on each other. In the opposite situation there is also a possibility to get result via \textbf{CE} results however, it will be necessary to change the algorithm and the dependence of how the functions depend on each other and add this to the algorithm when initiating optimization functions in the code. The same situation applies to  changing the objective function and the density used for approximation. Changing the objective function gives the possibility to use other methods. If it stays quadratic and the constraints are linear, quadratic programming techniques like the sequential quadratic programming method can still be used. If the objective function and the constraints stays convex, we use general methods from convex optimization. The \textbf{CE} method is based on well-known and simply classical rules and is therefore quite problem-independent. There are no rigid formulas and therefore it requires consideration for each problem separately. So it is very possible to reconsider it in another way, which can be interesting for further work. The \textbf{CE} method may also facilitate the analysis of modifications of epidemiological models related to virus mutation or delayed expression (v. \citeauthor{GubTayZhu2018:OCHMV}~\cite{GubTayZhu2018:OCHMV}, \citeauthor{KocGraLip2020:COVID}~\cite{KocGraLip2020:COVID}).
\vspace{6pt} 



\authorcontributions{KSz gives an idea of application \textbf{CE} to SIRC model; MKS realizes implementation and simulation.}

\funding{KSz are thankful to Faculty of Pure and Applied Mathematics for research funding.}


\conflictsofinterest{The authors declare no conflict of interest.}

\abbreviations{\label{abbrAll}The following abbreviations are used in this article:\\
\begin{tabular}{@{}lp{0.8\linewidth}}
\textbf{CE}  &  Cross Entropy Method (v. page \pageref{IS::important})\\
\textbf{FEM} & the finite element methods (v. page \pageref{FEM})\\
\textbf{GOF} & a goodness-of-fit (v. page \pageref{GOF})\\
\textbf{IS}  & the importance sampling (v. page \pageref{IS::important})\\
\textbf{SIR} & The \textbf{SIR} model is one of the simplest \emph{compartmental} models. The letters means the number of \textbf{S}usceptible-\textbf{I}nfected-\textbf{R}emoved  (and immune) or deceased individuals.  (v. \cite{KerMcK1991:MTE1}, \cite{HarLobMak2014:Exact}).\\
\textbf{SIRC}& The \textbf{SIR} model with  the additional group of partially resistant to the current strain people: \textbf{S}usceptible - \textbf{I}nfectious - \textbf{R}ecovered - \textbf{C}ross-immune  (v. page~\pageref{section:sircmodel}). 
\end{tabular}
}

\appendixtitles{yes} 
\appendix
\section{Optimization by the method of cross-entropy.} 
\subsection{\label{SumOfRank}Multiple selection to minimize the sum of ranks.}
Let $(a_{1},\dots,a_{N})$ be a permutation of integers $(1,2,\dots,N)$, where all permutations are equally probable. For every $i=1,2,\dots,N$ $Y_{i}$ will be the number of values $a_{1},\dots,a_{i}$, which are $\leqslant a_{i}$. $Y_{i}=\card\{1\leqslant j\leqslant i: a_j\leqslant i\}$ is the relative rank of the $i$-th object. The decision-maker observe the relative ranks and the realite ranks define filtration: $\mathcal{F}_i=\sigma\{Y_1,Y_2,\ldots,Y_i\}$. Let $\frak{S}$ be the set of Markov time with respect to $\{\mathcal{F}_i\}_{i=1}^N$. Define
	\begin{align*} 
	v&=\inf_{\tau} \bE (a_{\tau_{1}}+\dots+a_{\tau_{k}}), \qquad \text{where $\tau=(\tau_{1},\dots,\tau_{k})$, $\tau_i\in\frak{S}$.} 
	\end{align*}
	We want to find the optimal procedure $\tau^{*}=(\tau_{1}^{*},\dots,\tau_{k}^{*})$ and the value of the optimization problem $v$. 
	Let $\mathcal{F}_{(m)_{i}}$ be $\sigma$-algebra generated by $(Y_{1},\dots,Y_{m_{i}})$. If we accept
	\begin{align*} Z_{(m)_{k}}&=\bE (a_{\tau_{1}}+\dots+a_{\tau_{k}}|F_{(m)_{k}}), 
	\intertext{then}
	v&=\inf_{\tau}\bE Z_{\tau}, \qquad \tau=(\tau_{1},\dots,\tau_{k}). 
	\end{align*}
	The problem is reduced to the problem of stopping sequences multiple times $Z_{(m)_{k}}$.
	As shown in \cite{Nikolaev1} and \cite{Nikolaev2} the solution to the problem is the following strategy: \\
	If there are integer vectors
	\begin{equation}\label{wzor1}
	\begin{aligned}
	&\delta^{(k)}=(\delta^{(k)}_{1},\dots,\delta^{(k)}_{N-k+1}), 
	&0 \leqslant \delta^{(k)}_{1} \leqslant \dots \leqslant \delta^{(k)}_{N-k}  < \delta^{(k)}_{N-k+1} = N, \\
	&\delta^{(2)}=(\delta^{(2)}_{k-1},\dots,\delta^{(2)}_{N-1}),
	&0 \leqslant \delta^{(2)}_{k-1} \leqslant \dots \leqslant \delta^{(2)}_{N-2}  < \delta^{(2)}_{N-1} = N,\\
	&\delta^{(1)}=(\delta^{(1)}_{k},\dots,\delta^{(1)}_{N}),
	&0 \leqslant \delta^{(1)}_{k} \leqslant \dots \leqslant \delta^{(1)}_{N-1}  < \delta^{(1)}_{N} = N,\\
	&\delta^{i_{1}}_{j} \leqslant \delta^{i_{2}}_{j}, \qquad 1 \leqslant i_{1} \leqslant i_{2} \leqslant k, 
	& k-i_{1}+1\leqslant j \leqslant N-i_{2} + 1
	\end{aligned}
	\end{equation}
	then
	\[ 
	\begin{split}
	&\tau_{1}^{*}=min\{ m_{1}: y_{m_{1}} \leqslant \delta^{(k)}_{m_{1}} \}, \\
	&\tau_{i}^{*}=min\{ m_{i}>m_{i-1}: y_{m_{i}} \leqslant \delta^{(k-i+1)}_{m_{i}} \},
	\end{split} 
	\]
	on the set $F_{i-1}= \{ \omega: \tau_{1}^{*}=m_{1},\dots,\tau_{i-1}^{*}=m_{i-1}\}$. $i=2,\dots,k$, $F_{0}=\Omega $.
	For small $N$ there is a possibility to get accurate values by the analytical method. The \textbf{CE} method solves this problem by changing the estimation problem to the optimization problem.  Then by randomizing this problem using a defined family of probability density functions. With this the \textbf{CE} method solves this efficiently by making adaptive changes to this pdf and go in the direction of the theoretically optimal density.

\subsection{\label{MinMeanSumRank}Minimization of mean sum of ranks.} Let's consider the following problem of minimizing the mean sum of ranks:
	\begin{equation}\label{wzor2}
	\min_{x \in \chi} \bE S(x,R),
	\end{equation} 
	where $\chi=\{x=(x^{(1)},\dots,x^{(k)}):$ conditions from \eqref{wzor1} are preserved$\}$ is some defined set. $R=(R_{1},\dots,R_{N})$ is a random permutation of numbers $1,\dots,N$. $\hat{S}$ is an unbiased estimator $\bE S(x,R)$ with the following formula:
	\[
	\hat{S}(x)= \frac{1}{N_{1}}\sum_{n=1}^{N_{1}}(R_{n\tau_{1}}+\dots+R_{n\tau_{k}}),
	\]
	where $(R_{n1},\dots,R_{nN})$ is the nth copy of a random permutation $R$. \\
	Now a cross-entropy algorithm can be applied (see \ref{section:cemethod}, \cite{rubinstein1997optimization}). Let's define the indicator collections $\{I_{\{S(x)\leqslant \gamma \}} \}$ on $\chi$ for different levels $\gamma \in R$. Let $\{f((\cdotp; u)\}$ be the density family on $\chi$ parameterized by the actual parameter value $u$. For a specific $u$ we can combine \eqref{wzor2} with the problem of estimation
	\[ l(\gamma)=\bP_{u} (S(X)\leqslant \gamma)=\sum_{x}I_{\{S(x)\leqslant \gamma \}} f(x,u)=\bE_{u} I_{\{S(x)\leqslant \gamma \}}, \]
	where $\bP_{u}$ is a measure of the probability in which the random state $X$ has a density $\{f((\cdotp; u)\}$ and $\gamma$ is a known or unknown parameter. $l$ we estimate using the Kullback-Leibler distance.	The Kullback-Leibler distance is defined as follows:
    \[ 
    D(g,h)=\int g(x)\ln g(x) dx - \int g(x)h(x)dx 
    \]
    
     \noindent Usually $g(x)$ is chosen from the family of density $f(\cdotp;v)$. So here $D$ for $g$ and $f(x;v)$ it comes down to selecting such a reference parameter $v$ for which $-\int g(x)\ln f(x;v)dx$ is the smallest, that is, it comes down to maximization:
    \[ \max_{v} \int g(x) f(x;v) dx. \]
    After some transformations:
    \[ \max_{v}D(v)=\max_{v}\frac{1}{N}\sum_{i=1}^{N}\mathbb{I}_{\{S(X_{i})\geqslant \gamma \}}W(X_{i};u;w)\ln f(X_{i};v) \]
    
 \noindent So firstly, let's generate a pair $\{(\gamma_{t},u_{t})\}$, which we then update until the stopping criterion is met and the optimal pair is obtained $\{(\gamma^{*},u^{*})\}$. More precisely, arrange $u_{0}$ and choose not too small $\varrho$ and we proceed as follows:
\begin{description}
	\item[1) Updating $\gamma_{t}$]
	Generate a sample $X_{1},\dots,X_{N_{2}}$ from $\{f((\cdotp; u_{t-1})\}$. Calculate $\hat{S}(X_{1}),\dots,\hat{S}(X_{N_{2}})$ and sort in ascending order. For $\hat{\gamma}_{t}$ choose
	\[ \hat{\gamma}_{t}=\hat{S}_{(\lceil \varrho N \rceil)} \]
	\item[2) Updating $u_{t}$]
	$\hat{u}_{t}$ obtain from the Kullback-Leibler distance, that is, from maximization
	\begin{align}\label{wzor23}
	\max_{u}D(u)&=\max_{u}\bE_{u_{t-1}}I_{\{\hat{S}(x)\leqslant \gamma_{t} \}} ln f(X;u),\\[-2ex]
	\intertext{so}
	\label{wzor24}
	\max_{u}\hat{D}(u)&=\max_{u}\frac{1}{N_{2}}\sum_{n=1}^{N_{2}}I_{\{\hat{S}(x)\leqslant \hat{\gamma_{t}} \}} lnf(X_{n};u).
	\end{align}
	As in \cite{SofKroKeiNik2006:MBCP}, here a 3-dimensional matrix of parameters $u=\{u_{ijl}\}$ is consider.
	\begin{align*} 
	 u_{ijl}&=\Pr\{X_{j}^{(i)}=l\}, \text{  $i=1,\dots,k;$}\quad \text{ $j=k-i+1,\dots,N-i; l=0,\dots,N-1$.}\\[-2ex]
	\intertext{It seems that}
	f(x^{(i)}_{j};u)&=\sum_{l=0}^{N-1}u_{ijl} I_{\{x_{j}^{(i)}=l\}}.\\[-3ex]
	\intertext{and then after some transformations}
	\hat{u}^{(t)}_{ijl}&=\dfrac{\sum_{n=1}^{N_{2}}I_{\{\hat{S}(X_{n}) \leqslant \hat{\gamma}_{t} \} } W^{(t-1)}_{nij}I_{\{X_{nij}=l\} } }{\sum_{n=1}^{N_{2}}I_{\{\hat{S}(X_{n}) \leqslant \hat{\gamma}_{t} \} } }, \\
	W_{nij}^{(t-1)}&=\frac{\hat{u}_{ijX_{nij}}^{(0)} }{\hat{u}_{ijX_{nij} }^{(t-1)} },
	\end{align*}
	where $X_{n}=\{X_{nij}\}$, $X_{nij}$ is a random variable from $f(x_{j}^{(i)};\hat{u}_{t-1})$, corresponding to the formula \eqref{wzor24}.
	Instead of updating a parameter use the following smoothed version $$ \hat{u}_{t}=\alpha \hat{u}_{t} + (1-\alpha)\hat{u}_{t-1}.$$ 
	
	\item[3) Stopping Criterion]
    The criterion is from \cite{Rub1999:CEMethod}, which stop the algorithm when $\hat{\gamma_{T}}$ ($T$ is last step) has reached stationarity. To identify the stopping point of $T$, consider the following moving average process
    \[ B_{t}(K)=\frac{1}{K}\sum_{s=t-K+1}^{t}\hat{\gamma}_{t},\, t=s,s+1,\dots,\,s\geqslant K, \]
    where $K$ is fixed.
    \[ C_{t}(K)=\dfrac{\dfrac{1}{K-1}\{ \sum_{s=t-K+1}^{t} (\hat{\gamma}_{t}-B_{t}(K))^{2} \}}{B_{t}(K)^{2}}. \]
    Then let's define
    \[ C_{t}^{-}(K,R)=\min_{j=1,\dots,R}C_{t+j}(K) \]
    and
    \[ C_{t}^{+}(K,R)=\max_{j=1,\dots,R}C_{t+j}(K), \]
    where $R$ is fixed. \\
    Then the stopping criterion is defined as follows
    \begin{equation}
    T=\min\{ t:\dfrac{C_{t}^{+}(K,R)-C_{t}^{-}(K,R)}{C_{t}^{+}(K,R)}\leqslant \epsilon \}
    \end{equation}
    where $K$ and $R$ are fixed and $\epsilon$ is not too small.
\end{description}
		
\subsection{\label{VehRP} The vehicle routing problem.}
The classical vehicle routing problem (\textbf{VRP}) is defined on a graph $\mathbb{G} = (\mathbf{V}, \mathcal{A})$, where $\mathbf{V} = \{v_0, v_1,\ldots, v_n\}$ is a set of vertices and $\mathcal{A} = \{(v_i, v_j ): i, j \in \{0,\ldots, n\}$, $v_i, v_j \in \mathbf{V}\}$ is the arc set. A matrix $\mathbf{L} = (L_{i, j} )$ can be defined on $\mathcal{A}$, where the coefficient $L_{i, j}$ defines the distance between the nodes $v_i$ and $v_j$ and is proportional to the cost of travelling by the corresponding arc. There can be one or more vehicles starting off from the depot $v_0$ with a given capacity, visiting all or a subset of the vertices, and returning to the depot after having satisfied the demands at the vertices. The \emph{Stochastic Vehicle Routing Problem} (\textbf{SVRP}) arises when elements of the vehicle routing problem are stochastic – the set of customers visited, the demands at the vertices, or the travel times. 

Let us consider that a certain type of a product is distributed from a plant to $N$ customers, using a \emph{single} vehicle, having a fixed capacity $Q$. The vehicle strives to visit all the customers periodically to supply the product and replenish their inventories. On a given periodical trip through the network, on visiting a customer, an amount equal to the demand of that customer is downloaded from the vehicle, which then moves to the next site. The demands of a given customer during each period are modelled as independent and identically distributed random variables with known distribution. A reasonable assumption is that all the customers’ demands belong to a certain distribution (say normal) with varying parameters for different customers.

The predominant approach for solving the \textbf{SVRP} class of problems is to use a “here-and-now” optimization technique, where the sequence of customers to be visited is decided in advance. On the given route, if the vehicle fails to meet the demands of a customer, there is a \emph{recourse action} taken. The recourse action could be in the form of \emph{going back to the depot} to replenish and fulfil the customers' demand, and continue with the remaining sequence of the customers to be visited or any other meaningful formulation. The problem then reduces to a stochastic optimization problem where the sequence with the minimum expected distance of travel (or equivalently, the minimum expected cost) has to be arrived at. The alternative is to use a re-optimization strategy where upon failure at a node, the optimum route for the remaining nodes is recalculated. The degree of re-optimization varies. At one extreme is the use of a dynamic approach where one can re-optimize at any point, using the newly obtained data about customer demands, or to re-optimize after failure. Neuro Dynamic Programming has been used to implement techniques based on re-optimization (cf., e.g., \cite{Sec2000:Neuro}).

The vehicle departs at constant capacity Q and has no knowledge of the requirements it will encounter on the route, except for the probability distributions of individual orders. Hence, there is a positive probability that the vehicle runs out of the product along the route, in which case the remaining demands of that customer and the remaining customers further along the route are not satisfied (a failure route). Such failures are discouraged with penalties, which are functions of the recourse actions taken. Each customer in the set can have a unique penalty cost for not satisfying the demand. The cost function for a particular route travelled by the vehicle during a period is calculated as the sum of all the arcs visited and the penalties (if any) imposed. If the vehicle satisfies all the demands on that route, the cost of that route will simply be the sum of the arcs visited including the arc from the plant to the first customer visited, and the arc from the last customer visited back to the plant. 

Alternatively, if the vehicle fails to meet the demands of a particular customer, the vehicle heads back to the plant at that point, terminating the remaining route. The cost function is then the sum of all the arcs visited (including the arc from the customer where the failure occurred back to the plant) and the penalty for that customer. In addition, the penalties for the remaining customers who were not visited will also be imposed. Thus, a given route can have a range of cost function values associated with it. The objective is to find the route for which the expected value of the cost function is minimum compared to all other routes.

Let us adopt the “here-and-now” optimization approach. It means that the operator decides the sequence of vertices to be visited in advance and independent of the demands encountered. Such approach leads to a \emph{discrete stochastic optimization problem} with respect to the discrete set of finite routes that can be taken. Let $G(r) := \mathbf{E}[H(r, D_r )]$, where $H(r,D_r)$ is the deterministic cost function of the route $r=(r_1,\ldots,r_N)$, $R$ -- is the discrete and finite (or countable finite) feasible set of values that r can take, and demands $D_{r_i}$ are a random variable that may or may not depend on the parameters $r_i$. The class of discrete stochastic optimization problems we are considering here are those of the form
\begin{align*}
\min_{r\in R}\{G(r)\}.
\end{align*}
Thus, $H(r, D_r)$ is a random variable whose expected value, $G(r)$ is usually estimated by Monte Carlo simulation. The global optimum solution set can then be denoted by 
\begin{align*}
R^* = \{r*\in R : G(r^*) \leq G(r ), \forall_{ r\in R}\}.
\end{align*}
The determination of $R^*$ can be done by \textbf{CE} (v. Rubinstein \cite{rubinstein1997optimization,Rub1999:CEMethod,Rub2002:CERE}). In this case also the basic idea is to connect the underlying optimization problem to a problem of the estimating rare-event probabilities (v. \cite{de2005tutorial}  in the context of the deterministic optimization). Here the \textbf{CE} method is used in the context of the \emph{discrete stochastic optimization} and the Monte Carlo techniques are needed to estimate the objective function. In solving the problem by \textbf{CE} methods some practical issues should be solved, as to when to draw new samples and how many samples to use.



\reftitle{References}



\bibliography{\jobname}

\begin{thebibliography}{-------}
\providecommand{\natexlab}[1]{#1}

\bibitem[Martino \em{et~al.}(2018)Martino, Luengo, and M{\'i}guez]{Martino2018}
Martino, L.; Luengo, D.; M{\'i}guez, J., Introduction.
\newblock In {\em Independent Random Sampling Methods}; Springer International
  Publishing: Cham,  2018; pp. 1--26.
\newblock
  doi:{\changeurlcolor{black}\href{https://doi.org/10.1007/978-3-319-72634-2_1}{\detokenize{10.1007/978-3-319-72634-2_1}}}.

\bibitem[Metropolis(1987)]{Met1987:BeginMC}
Metropolis, N.
\newblock The beginning of the Monte Carlo method.
\newblock {\em Los Alamos Sci.} {\bf 1987}, {\em 15},~125--130.
\newblock Special Issue: Stanislaw Ulam 1909–1984.

\bibitem[{Metropolis} and {Ulam}(1949)]{MetUla1949:MCM}
{Metropolis}, N.; {Ulam}, S.
\newblock {The Monte Carlo method}.
\newblock {\em {J. Am. Stat. Assoc.}} {\bf 1949}, {\em 44},~335--341.
\newblock \ZBL{0033.28807},
  doi:{\changeurlcolor{black}\href{https://doi.org/10.2307/2280232}{\detokenize{10.2307/2280232}}}.

\bibitem[Marshall(1956)]{Mar1956:IS}
Marshall, A.
\newblock The use of multistage sampling schemes in Monte Carlo computations.
\newblock  Symposium on Monte Carlo; Wiley: New York,  1956; p. 123–140.

\bibitem[Rubinstein(1997)]{rubinstein1997optimization}
Rubinstein, R.Y.
\newblock Optimization of computer simulation models with rare events.
\newblock {\em European Journal of Operational Research} {\bf 1997}, {\em
  99},~89--112.
\newblock
  doi:{\changeurlcolor{black}\href{https://doi.org/10.1016/S0377-2217(96)00385-2}{\detokenize{10.1016/S0377-2217(96)00385-2}}}.

\bibitem[Sani and Kroese(2007)]{sani2007optimal}
Sani, A.; Kroese, D.
\newblock {Optimal Epidemic Intervention of {HIV} Spread Using Cross-Entropy
  Method}.
\newblock  Proceedings of the International Congress on Modelling and
  Simulation (MODSIM); Oxley, L.; Kulasiri, D., Eds.; Modelling and Simulation
  Society of Australia and New Zeeland: Christchurch, New Zeeland,  2007; pp.
  448--454.

\bibitem[Asamoah \em{et~al.}(2018)Asamoah, Nyabadza, Seidu, Chand, and
  Dutta]{AsoNya2018:Meningitis}
Asamoah, J.; Nyabadza, F.; Seidu, B.; Chand, M.; Dutta, H.
\newblock Mathematical Modelling of Bacterial Meningitis Transmission Dynamics
  with Control Measures.
\newblock {\em Computational and Mathematical Methods in Medicine} {\bf 2018},
  {\em Article ID 2657461},~21 page.
\newblock \PMID{29780431}, \PMCID{PMC5892307},
  doi:{\changeurlcolor{black}\href{https://doi.org/10.1155/2018/2657461}{\detokenize{10.1155/2018/2657461}}}.

\bibitem[Vereen(2008)]{Ver2008:SIRC}
Vereen, K.
\newblock An SCIR Model of Meningococcal Meningitis.
\newblock Master's thesis, Virginia Commonwealth University, Richmond, Virginia
  23284,  2008.
\newblock
  \href{https://scholarscompass.vcu.edu/cgi/viewcontent.cgi?article=1709&context=etd}{Electronic
  theses and dissertations (ETDs) of VCU.}

\bibitem[Casagrandi \em{et~al.}(2006)Casagrandi, Bolzoni, Levin, and
  Andreasen]{casagrandi2006sirc}
Casagrandi, R.; Bolzoni, L.; Levin, S.A.; Andreasen, V.
\newblock {The {SIRC} model and influenza {A}}.
\newblock {\em Math. Biosci.} {\bf 2006}, {\em 200},~152--169.
\newblock \MR{2225509}, \PMID{16504214},
  doi:{\changeurlcolor{black}\href{https://doi.org/10.1016/j.mbs.2005.12.029}{\detokenize{10.1016/j.mbs.2005.12.029}}}.

\bibitem[{Parry}(1969)]{Par1969:Entropy}
{Parry}, W.
\newblock {\em {Entropy and generators in ergodic theory}}; Mathematics Lecture
  Note Series, W. A. Benjamin, Inc.: New York-Amsterdam,  1969; pp. xii, 124 p.
\newblock \ZBL{0175.34001}.

\bibitem[Cover and Thomas(2006)]{CovTho2006:ElementIT}
Cover, T.M.; Thomas, J.A.
\newblock {\em Elements of information theory}, second ed.; Wiley-Interscience
  [John Wiley \& Sons], Hoboken, NJ,  2006; pp. xxiv+748.
\newblock
  doi:{\changeurlcolor{black}\href{https://doi.org/10.1002/047174882X}{\detokenize{10.1002/047174882X}}}.

\bibitem[{Kullback} and {Leibler}(1951)]{KulLei1951:Information}
{Kullback}, S.; {Leibler}, R.A.
\newblock {On information and sufficiency}.
\newblock {\em {Ann. Math. Stat.}} {\bf 1951}, {\em 22},~79--86.
\newblock \ZBL{0042.38403},
  doi:{\changeurlcolor{black}\href{https://doi.org/10.1214/aoms/1177729694}{\detokenize{10.1214/aoms/1177729694}}}.

\bibitem[Inglot(2014)]{Ing2014:TIandSM}
Inglot, T.
\newblock Teoria informacji a statystyka matematyczna.
\newblock {\em Mathematica Applicanda} {\bf 2014}, {\em 42},~115--174.
\newblock
  doi:{\changeurlcolor{black}\href{https://doi.org/10.14708/ma.v42i1.521}{\detokenize{10.14708/ma.v42i1.521}}}.

\bibitem[{Csisz\'ar}(1967)]{Csi1967:Information}
{Csisz\'ar}, I.
\newblock {Information-type measures of difference of probability distributions
  and indirect observations}.
\newblock {\em {Stud. Sci. Math. Hung.}} {\bf 1967}, {\em 2},~299--318.
\newblock \ZBL{0157.25802}.

\bibitem[Amari(1985)]{Ama1985:Diff}
Amari, S.i.
\newblock {\em Differential-geometrical methods in statistics}; Vol.~28, {\em
  Lecture Notes in Statistics}, Springer-Verlag, New York,  1985; pp. v+290.
\newblock \MR{788689},
  doi:{\changeurlcolor{black}\href{https://doi.org/10.1007/978-1-4612-5056-2}{\detokenize{10.1007/978-1-4612-5056-2}}}.

\bibitem[{Rubinstein}(1999)]{Rub1999:CEMethod}
{Rubinstein}, R.
\newblock The cross-entropy method for combinatorial and continuous
  optimization.
\newblock {\em Methodol. Comput. Appl. Probab.} {\bf 1999}, {\em 1},~127--190.
\newblock \ZBL{0941.65061},
  doi:{\changeurlcolor{black}\href{https://doi.org/10.1023/A:1010091220143}{\detokenize{10.1023/A:1010091220143}}}.

\bibitem[{ {Ferguson}, T.S.}(1989)]{Fer1989:Who}
{ {Ferguson}, T.S.}.
\newblock Who solved the secretary problem?
\newblock {\em Statistical Science} {\bf 1989}, {\em 4},~282--296.

\bibitem[{ {Szajowski}, K.}(1982)]{sza1982:Ath}
{ {Szajowski}, K.}.
\newblock Optimal choice problem of a-th object.
\newblock {\em {Matem. Stos.}} {\bf 1982}, {\em 19},~51--65.
\newblock in Polish,
  doi:{\changeurlcolor{black}\href{https://doi.org/10.14708/ma.v10i19.1533}{\detokenize{10.14708/ma.v10i19.1533}}}.

\bibitem[Polushina(2010)]{Pol2010:CEBCP}
Polushina, T.V.
\newblock Estimating Optimal Stopping Rules in the Multiple Best Choice Problem
  with Minimal Summarized Rank via the Cross-Entropy Method. In {\em
  Exploitation of Linkage Learning in Evolutionary Algorithms}; {Chen,
  Ying-Ping}., Ed.; Springer Berlin Heidelberg: Berlin, Heidelberg,  2010; pp.
  227--241.
\newblock
  doi:{\changeurlcolor{black}\href{https://doi.org/10.1007/978-3-642-12834-9_11}{\detokenize{10.1007/978-3-642-12834-9_11}}}.

\bibitem[Stachowiak(2019)]{Sta2019}
Stachowiak, M.
\newblock The cross-entropy method and its applications.
\newblock Technical report, Faculty of Pure and Applied Mathematics, Wrocław
  University of Science and Technology, Wroclaw,  2019.
\newblock 33p. Engineering Thesis.

\bibitem[Dror(2002)]{Dro2002:VRSD}
Dror, M.
\newblock Vehicle Routing with Stochastic Demands: Models \& Computational
  Methods. In {\em Modeling Uncertainty}; {Dror, M. and {L’Ecuyer}, P. and
  Szidarovszky, F.}., Ed.; Springer: New York, NY,  2002; Vol.~46, {\em
  International Series in Operations Research \& Management Science}, pp.
  625--649.
\newblock
  doi:{\changeurlcolor{black}\href{https://doi.org/10.1007/0-306-48102-2_25}{\detokenize{10.1007/0-306-48102-2_25}}}.

\bibitem[Chepuri and Homem-de Mello(2005)]{CheHom2005:vehicle}
Chepuri, K.; Homem-de Mello, T.
\newblock Solving the vehicle routing problem with stochastic demands using the
  cross-entropy method.
\newblock {\em Ann. Oper. Res.} {\bf 2005}, {\em 134},~153--181.
\newblock
  doi:{\changeurlcolor{black}\href{https://doi.org/10.1007/s10479-005-5729-7}{\detokenize{10.1007/s10479-005-5729-7}}}.

\bibitem[Ekeland and Temam(1999)]{ekeland1999convex}
Ekeland, I.; Temam, R.
\newblock {\em Convex analysis and variational problems}; Vol.~28, SIAM,  1999.

\bibitem[Glowinski(2008)]{glowinski2008lectures}
Glowinski, R.
\newblock {\em Lectures on numerical methods for non-linear variational
  problems}; Springer Science and Business Media,  2008.

\bibitem[Mumford and Shah(1989)]{mumford1989optimal}
Mumford, D.; Shah, J.
\newblock Optimal approximations by piecewise smooth functions and associated
  variational problems.
\newblock {\em Communications on Pure and Applied Mathematics} {\bf 1989}, {\em
  42},~577--685.

\bibitem[Klimov \em{et~al.}(1999)Klimov, Simonsen, Fukuda, and
  Cox]{klimov1999surveillance}
Klimov, A.; Simonsen, L.; Fukuda, K.; Cox, N.
\newblock Surveillance and impact of influenza in the United States.
\newblock {\em Vaccine} {\bf 1999}, {\em 17},~S42--S46.
\newblock Suplement 1. \PMID{10471179},
  doi:{\changeurlcolor{black}\href{https://doi.org/10.1016/S0264-410X(99)00104-8}{\detokenize{10.1016/S0264-410X(99)00104-8}}}.

\bibitem[Simonsen \em{et~al.}(1997)Simonsen, Clarke, Williamson, Stroup, Arden,
  and Schonberger]{simonsen1997impact}
Simonsen, L.; Clarke, M.J.; Williamson, G.D.; Stroup, D.F.; Arden, N.H.;
  Schonberger, L.B.
\newblock The impact of influenza epidemics on mortality: introducing a
  severity index.
\newblock {\em American Journal of Public Health} {\bf 1997}, {\em
  87},~1944--1950.
\newblock \PMID{9431281},
  doi:{\changeurlcolor{black}\href{https://doi.org/10.2105/AJPH.87.12.1944}{\detokenize{10.2105/AJPH.87.12.1944}}}.

\bibitem[Earn \em{et~al.}(2002)Earn, Dushoff, and Levin]{earn2002ecology}
Earn, D.J.; Dushoff, J.; Levin, S.A.
\newblock Ecology and evolution of the flu.
\newblock {\em Trends in Ecology \& Evolution} {\bf 2002}, {\em 17},~334--340.

\bibitem[Andreasen \em{et~al.}(1997)Andreasen, Lin, and
  Levin]{AndLinLev1997:Dynamics}
Andreasen, V.; Lin, J.; Levin, S.A.
\newblock The dynamics of cocirculating influenza strains conferring partial
  cross-immunity.
\newblock {\em J. Math. Biol.} {\bf 1997}, {\em 35},~825--842.
\newblock \MR{1479341},
  doi:{\changeurlcolor{black}\href{https://doi.org/10.1007/s002850050079}{\detokenize{10.1007/s002850050079}}}.

\bibitem[Lin \em{et~al.}(1999)Lin, Andreasen, and Levin]{lin1999dynamics}
Lin, J.; Andreasen, V.; Levin, S.A.
\newblock {Dynamics of influenza {A} drift: the linear three-strain model}.
\newblock {\em Mathematical Biosciences} {\bf 1999}, {\em 162},~33--51.

\bibitem[Iacoviello and Stasio(2013)]{iacoviello2013optimal}
Iacoviello, D.; Stasio, N.
\newblock {Optimal control for {SIRC} epidemic outbreak}.
\newblock {\em Computer Methods and Programs in Biomedicine} {\bf 2013}, {\em
  110},~333--342.
\newblock \PMID{23399104},
  doi:{\changeurlcolor{black}\href{https://doi.org/10.1016/j.cmpb.2013.01.006}{\detokenize{10.1016/j.cmpb.2013.01.006}}}.

\bibitem[Zaman \em{et~al.}(2008)Zaman, Kang, and Jung]{zaman2008stability}
Zaman, G.; Kang, Y.H.; Jung, I.H.
\newblock {Stability analysis and optimal vaccination of an {SIR} epidemic
  model}.
\newblock {\em BioSystems} {\bf 2008}, {\em 93},~240--249.

\bibitem[Kamien and Schwartz(1991)]{kamien2012dynamic}
Kamien, M.I.; Schwartz, N.L.
\newblock {\em {Dynamic Optimization: {T}he Calculus of Variations and Optimal
  Control in Economics and Management}}, second ed.; Vol.~31, {\em Advanced
  Textbooks in Economics}, North-Holland Publishing Co., Amsterdam,  1991; pp.
  xviii+377.
\newblock \MR{1159711}.

\bibitem[Fang \em{et~al.}(2002)Fang, Tsuchiya, and Yamamoto]{fang2002finite}
Fang, Q.; Tsuchiya, T.; Yamamoto, T.
\newblock Finite difference, finite element and finite volume methods applied
  to two-point boundary value problems.
\newblock {\em J. Comput. Appl. Math.} {\bf 2002}, {\em 139},~9--19.
\newblock \MR{1876870},
  doi:{\changeurlcolor{black}\href{https://doi.org/10.1016/S0377-0427(01)00392-2}{\detokenize{10.1016/S0377-0427(01)00392-2}}}.

\bibitem[Caglar \em{et~al.}(2006)Caglar, Caglar, and Elfaituri]{caglar2006}
Caglar, H.; Caglar, N.; Elfaituri, K.
\newblock B-spline interpolation compared with finite difference, finite
  element and finite volume methods which applied to two-point boundary value
  problems.
\newblock {\em Applied Mathematics and Computation} {\bf 2006}, {\em
  175},~72--79.
\newblock
  doi:{\changeurlcolor{black}\href{https://doi.org/10.1016/j.amc.2005.07.019}{\detokenize{10.1016/j.amc.2005.07.019}}}.

\bibitem[Papadopoulos and Giovanis(2018)]{PapGio2018:SFEM}
Papadopoulos, V.; Giovanis, D.G.
\newblock {\em Stochastic finite element methods}; Mathematical Engineering,
  Springer, Cham,  2018; pp. xxi+138.
\newblock An introduction,
  doi:{\changeurlcolor{black}\href{https://doi.org/10.1007/978-3-319-64528-5}{\detokenize{10.1007/978-3-319-64528-5}}}.

\bibitem[Lenhart and Workman(2007)]{lenhart2007optimal}
Lenhart, S.; Workman, J.T.
\newblock {\em Optimal control applied to biological models}; CRC Press,  2007.

\bibitem[Hazelbag \em{et~al.}(2020)Hazelbag, Dushoff, Dominic, Mthombothi, and
  Delva]{hazelbag2020calibration}
Hazelbag, C.M.; Dushoff, J.; Dominic, E.M.; Mthombothi, Z.E.; Delva, W.
\newblock Calibration of individual-based models to epidemiological data: A
  systematic review.
\newblock {\em PLoS Computational Biology} {\bf 2020}, {\em 16},~e1007893.

\bibitem[Taylor \em{et~al.}(2010)Taylor, Pawar, Kruzikas, Gilmore, Pandya,
  Iskandar, and Weinstein]{taylor2010methods}
Taylor, D.C.; Pawar, V.; Kruzikas, D.; Gilmore, K.E.; Pandya, A.; Iskandar, R.;
  Weinstein, M.C.
\newblock Methods of model calibration.
\newblock {\em Pharmacoeconomics} {\bf 2010}, {\em 28},~995--1000.

\bibitem[Taynitskiy \em{et~al.}(2017)Taynitskiy, Gubar, and
  Zhu]{TayGubZhu2017:OICESHM}
Taynitskiy, V.; Gubar, E.; Zhu, Q.
\newblock Optimal Impulsive Control of Epidemic Spreading of Heterogeneous
  Malware.
\newblock {\em IFAC-PapersOnLine} {\bf 2017}, {\em 50},~15038 -- 15043.
\newblock 20th IFAC World Congress,
  doi:{\changeurlcolor{black}\href{https://doi.org/10.1016/j.ifacol.2017.08.2515}{\detokenize{10.1016/j.ifacol.2017.08.2515}}}.

\bibitem[Gubar \em{et~al.}(2018)Gubar, Taynitskiy, and
  Zhu]{GubTayZhu2018:OCHMV}
Gubar, E.; Taynitskiy, V.; Zhu, Q.
\newblock Optimal Control of Heterogeneous Mutating Viruses.
\newblock {\em Games} {\bf 2018}, {\em 9}.
\newblock
  doi:{\changeurlcolor{black}\href{https://doi.org/10.3390/g9040103}{\detokenize{10.3390/g9040103}}}.

\bibitem[Kochańczyk \em{et~al.}(2020)Kochańczyk, Grabowski, and
  Lipniacki]{KocGraLip2020:COVID}
Kochańczyk, M.; Grabowski, F.; Lipniacki, T.
\newblock Dynamics of COVID-19 pandemic at constant and time-dependent contact
  rates.
\newblock {\em Math. Model. Nat. Phenom.} {\bf 2020}, {\em 15},~28.
\newblock
  doi:{\changeurlcolor{black}\href{https://doi.org/10.1051/mmnp/2020011}{\detokenize{10.1051/mmnp/2020011}}}.

\bibitem[Kermack and McKendrick(1991)]{KerMcK1991:MTE1}
Kermack, W.; McKendrick, A.
\newblock Contributions to the mathematical theory of epidemics--I. 1927.
\newblock {\em Bull Math Biol.} {\bf 1991}, {\em 53},~35--55.
\newblock Reprint of the Proc. of the Royal Society of London. Series A, 1927,
  Vol. 115, No. 772, pp. 700-721. \doi{10.1098/rspa.1927.0118}. The same
  journal published the parts II and III of the article in 1932 and 1933,
  respectively. \PMID{2059741}; \JFM{53.0517.01},
  doi:{\changeurlcolor{black}\href{https://doi.org/10.1007/BF02464423}{\detokenize{10.1007/BF02464423}}}.

\bibitem[Harko \em{et~al.}(2014)Harko, Lobo, and Mak]{HarLobMak2014:Exact}
Harko, T.; Lobo, F.S.N.; Mak, M.K.
\newblock Exact analytical solutions of the
  {Susceptible}-{Infected}-{Recovered} ({SIR}) epidemic model and of the {SIR}
  model with equal death and birth rates.
\newblock {\em Applied Mathematics and Computation} {\bf 2014}, {\em 236},~184
  -- 194.
\newblock \MR{3197716},
  doi:{\changeurlcolor{black}\href{https://doi.org/10.1016/j.amc.2014.03.030}{\detokenize{10.1016/j.amc.2014.03.030}}}.

\bibitem[{{Nikolaev}, M. L. }(1998)]{Nikolaev1}
{{Nikolaev}, M. L. }.
\newblock {On optimal multiple stopping of Markov sequences.}
\newblock {\em Theory Probab. Appl.} {\bf 1998}, {\em 43},~298--306.
\newblock Translation from Teor. Veroyatn. Primen. 43, No. 2, 374-382 (1998).
  \ZBL{0971.60046},
  doi:{\changeurlcolor{black}\href{https://doi.org/10.1137/S0040585X9797691X}{\detokenize{10.1137/S0040585X9797691X}}}.

\bibitem[{Nikolaev}(1998)]{Nikolaev2}
{Nikolaev}, M.L.
\newblock {Optimal multi-stopping rules}.
\newblock {\em {Obozr. Prikl. Prom. Mat.}} {\bf 1998}, {\em 5},~309--348.
\newblock \ZBL{1075.91511}.

\bibitem[{Safronov} \em{et~al.}(2006){Safronov}, {Kroese}, {Keith}, and
  {Nikolaev}]{SofKroKeiNik2006:MBCP}
{Safronov}, G.Y.; {Kroese}, D.P.; {Keith}, J.M.; {Nikolaev}, M.L.
\newblock {Simulations of thresholds in multiple best choice problem}.
\newblock {\em {Obozr. Prikl. Prom. Mat.}} {\bf 2006}, {\em 13},~975--982.
\newblock \ZBL{1161.60333}.

\bibitem[Secomandi(2000)]{Sec2000:Neuro}
Secomandi, N.
\newblock Comparing neuro-dynamic programming algorithms for the vehicle
  routing problem with stochastic demands.
\newblock {\em Computers \& Operations Research} {\bf 2000}, {\em 27},~1201 --
  1225.
\newblock
  doi:{\changeurlcolor{black}\href{https://doi.org/10.1016/S0305-0548(99)00146-X}{\detokenize{10.1016/S0305-0548(99)00146-X}}}.

\bibitem[{Rubinstein}(2002)]{Rub2002:CERE}
{Rubinstein}, R.Y.
\newblock {Cross-entropy and rare events for maximal cut and partition
  problems}.
\newblock {\em {ACM Trans. Model. Comput. Simul.}} {\bf 2002}, {\em
  12},~27--53.
\newblock \ZBL{1390.90482},
  doi:{\changeurlcolor{black}\href{https://doi.org/10.1145/511442.511444}{\detokenize{10.1145/511442.511444}}}.

\bibitem[de~Boer \em{et~al.}(2005)de~Boer, Kroese, Mannor, and
  Rubinstein]{de2005tutorial}
de~Boer, P.T.; Kroese, D.P.; Mannor, S.; Rubinstein, R.Y.
\newblock A tutorial on the cross-entropy method.
\newblock {\em Ann. Oper. Res.} {\bf 2005}, {\em 134},~19--67.
\newblock \MR{2136658},
  doi:{\changeurlcolor{black}\href{https://doi.org/10.1007/s10479-005-5724-z}{\detokenize{10.1007/s10479-005-5724-z}}}.

\end{thebibliography}
\sampleavailability{The implementation codes for the algorithms used in the examples are available from the authors.}


\end{document}